\documentclass[preprint,12pt]{elsarticle}

\usepackage{xurl}
\usepackage{amssymb}
\usepackage{amsmath}
\usepackage{amsthm}
\usepackage{algorithm}
\usepackage{algpseudocode}
\usepackage{booktabs}
\usepackage{graphicx}
\usepackage{bm}
\usepackage{xcolor}
\usepackage[normalem]{ulem}
\usepackage[capitalize]{cleveref}
\usepackage{epstopdf}
\usepackage[section]{placeins}
\usepackage{tikz}
\usepackage{subcaption}
\usepackage{siunitx}

\usepackage{comment}

\graphicspath{{Figs/}{./}}

\newcommand{\csrc}{Beijing Computational Science Research Center, Beijing, China}
\newcommand{\odu}{Department of Mathematics \& Statistics,
  Old Dominion University,\\
  Norfolk, VA 23508, USA}

\usepackage{unicode-math-input}
\let\cev\overleftarrow

\newcommand{\Rey}{\mathrm{Re}}
\newcommand{\May}{\mathrm{Ma}}
\newcommand{\Kn}{\mathrm{Kn}}
\newcommand{\cs}{c_{\mathrm{s}}}
\newcommand{\diff}{\mathrm{d}}
\newcommand{\ssstyle}{\scriptscriptstyle}
\newcommand{\fazero}{f_\alpha^{\ssstyle (0)}}
\newcommand{\faone}{f_\alpha^{\ssstyle (1)}}
\newcommand{\fatwo}{f_\alpha^{\ssstyle (2)}}

\newcommand{\tensorfont}[1]{\bm{\mathsf{#1}}}
\newcommand{\tensorgreek}[1]{\pmb{#1}}

\newcommand{\Adelta}{\tensorfont{A}_{\delta}}
\newcommand{\Tdelta}{\tensorfont{T}_{\!\!\delta}}

\newcommand{\qvector}[1]{\bm{\mathrm{#1}}}

\newcommand{\CD}{C_{\ssstyle\mathrm{D}}}
\newcommand{\CL}{C_{\ssstyle\mathrm{L}}}
\newcommand{\FD}{F_{\ssstyle\mathrm{D}}}
\newcommand{\FL}{F_{\ssstyle\mathrm{L}}}

\newcommand{\fceva}{f_{\cev{\alpha}}}
\newcommand{\Omegaceva}{\Omega_{\cev{\alpha}}}

\begin{document}

\begin{frontmatter}

\title{Residual-Driven Lattice Boltzmann Equation. 
\\Part I: Framework and Nonlinear Solution Methods}

\author[1]{Xingzhou Ji}
\ead{jixz@csrc.ac.cn}
\author[2]{Li-Shi Luo\corref{cor1}}
\ead{lluo@odu.edu}

\affiliation[1]{organization={\csrc}}
\affiliation[2]{organization={\odu}}
\cortext[cor1]{Corresponding author}

\begin{abstract}
We developed a residual-driven framework of the lattice Boltzmann equation (LBE)
to integrate general nonlinear iterative methods. The nonlinear residual is constructed directly from the complete LBE collision-streaming update, allowing general nonlinear solvers to be integrated without modifying the original LBE algorithmic structure, \emph{i.e.}, the kernel consisting of collision-streaming, boundary conditions, and various forcing and interaction effects.
Sophisticated techniques, such as weighted fixed-point iteration, 
full approximation scheme (FAS) multigrid, and Jacobian-free Newton--Krylov methods 
are integrated into the LBE via this framework,
and preconditioned LBE is further incorporated as a nonlinear smoother to
demonstrate solver composability.

Numerical tests validate that FAS and Newton methods do converge to the same discrete
steady fixed-point obtained with the primitive LBE iteration. For the 2D lid-driven cavity flow, 
the FAS scheme reduces the measured solution time by a factor of approximately 44 at
$\Rey=1000$ and more than 224 at $\Rey=5000$. Accurate results are also obtained
for the steady flows past a cylinder at $\Rey=20$ and the NACA0012 airfoil at $\Rey=1000$. 
The low Mach-number tests for the flow past a NACA0012 airfoil further show that
both the FAS scheme and preconditioning can not only substantially reduce the work 
required to reach steady solutions, but also mitigate the zero Mach number crisis.
\end{abstract}


\end{frontmatter}

\section{\label{sec:1} Introduction}
The lattice Boltzmann equation (LBE) has evolved into a widely recognized numerical method for simulating nearly incompressible flows (cf., \emph{e.g.}, the review \cite{lallemand2021lattice} and refs. therein). Its main advantages include a simple algorithmic structure and straightforward implementation on modern high-performance computing (HPC) architectures. This characteristic of the LBE has led to its broad adoption in areas of computational fluid dynamics (CFD) where direct numerical simulations (DNS)
are desirable in the cases such as 
turbulent flows and/or complex fluids in which fine spatiotemporal resolutions are required;
this is also rooted in its low numerical dissipation feature which is critical for DNS applications. 
The efficiency of the LBE algorithm is, at least, in part due to its
relatively light in floating-point operations (FLOPs) per iteration.
One striking attribute of the LBE is its \emph{coherent} discretization of phase space 
$\Gamma := (\bm{r},\, \bm{\xi})$ and time $t$ which seems to rigidly bind its time 
step size $\delta_t$ with its spatial resolution $\delta_x$ together
\cite{HeXY1997prer,HeXY1997pre}.
Consequently, the LBE may be seen inherently inefficient for steady-state problems or 
long-time integrations, which prefer adaptive time stepping and other more sophisticated 
time integration techniques.
To appreciate this point, let us consider the dimensionless 
discrete velocity Boltzmann equation (DVBE) in the incompressible 
limit \cite{bardos1991jsp}, from which the LBE may also be derived:
\begin{equation}
	\frac{\partial f_\alpha}{\partial t}
	+ \frac{1}{\epsilon}\bm{\xi}_\alpha \cdot\bm{\nabla} f_\alpha
	=
	\frac{1}{\epsilon^2}\Omega_\alpha,
\quad 
	\epsilon := \Kn = \frac{\May}{\Rey}\rightarrow 0 ,
\label{eqn:dvbe0}
\end{equation}
where $f_\alpha(\bm{r},\, t)$ is the particle velocity distribution function associated with the discrete velocity $\bm{\xi}_\alpha$, 
and $\Omega_\alpha$ is the collision term. Under diffusive scaling, both the transport and collision terms 
become stiff \cite{Inamuro1997pof, junk2005asymptotic}.
A conventional implicit-explicit (IMEX) splitting with explicit transport may therefore remain restrictive \cite{dimarco2027cicp}. In general, the primitive LBE with a uniform spatial resolution does not appear to be particularly well suited for directly solving \cref{eqn:dvbe0} in this stiff limit.

There have been persistent, though sporadic, efforts to integrate advanced time-stepping techniques with the LBE.
The multigrid (MG) methods have been introduced into the LBE to accelerate convergence 
for steady states \cite{tolke1998implicit,mavriplis2006multigrid}, 
and the dual-time-stepping (DTS) strategy \cite{gsell2020multigrid}
has extended this idea to transient flows. 
In addition, the Newton--Krylov--Schwarz methods \cite{huang2015fully}, the preconditioned inexact Newton methods \cite{huang2016nonlinearly}, and the MG-preconditioned Newton methods \cite{tolke1998implicit} have 
all been incorporated with LBE to 
substantially enhance its computational efficiency.
These studies demonstrated that the LBE iteration can be effectively embedded in more general solution procedures. 
However, most existing approaches \cite{tolke1998implicit,huang2015fully,huang2016nonlinearly}
constructed the nonlinear residuals based on the DVBE \eqref{eqn:dvbe0}, as opposed to the LBE \emph{per se},
thereby altering the discretizations defined by the LBE, leading to different 
finite-difference (FD) \cite{tolke1998implicit,huang2015fully} 
or finite-volume (FV) schemes \cite{li2019conservation}, 
and consequently forfeiting some, if not all, known desirable characteristics 
of the original LBE algorithm.

The importance of delineating the origin of discretizations 
so to distinguish and maintain those derived from the LBE
cannot be overstated.
The LBE has established
a broad range of physics models based on its versatile collision modeling
(cf., \emph{e.g.}, \cite{Ladd1993prl, Dunweg2007pre, Dellar2002jcp, Dellar2014siamjsc, lallemand2021lattice}) and 
its capability for large-scale 
and challenging industrial CFD applications 
(cf. \emph{e.g.}, \cite{konig2018exa,kiris2023hlpw}),
which may all be credited to its advantageous discretizations. The discretized DVBE differs from the LBE in terms of their discrete operators
and thus their numerical characteristics. Simply put, the advantages of the LBE, such as 
its low dissipation, isotropy, and efficiency (cf., \emph{e.g.}, \cite{Lallemand2000pre}), 
not to mention its implementation simplicity, are directly akin to its algorithm core --- the collision-streaming structure. 
Therefore, the objective of this work is to establish a framework to enable the LBE to 
integrate \emph{any} advanced nonlinear solution techniques without affecting its 
algorithmic kernel encompassing collision-streaming, boundary conditions, and various 
interactions and forcing effects. Sophisticated nonlinear iteration techniques will only 
serve as a wrap about the existing LBE algorithmic core, thus enhancing the LBE efficiency 
without compromising or sacrificing its advantages.
To attain this objective, the entire LBE algorithmic core must be encapsulated 
and preserved as a
whole entity. The evolution of the LBE kernel serves two purposes: it provides (a) the 
residual computed from the distribution functions to an outer nonlinear solver or
a time-stepping algorithm wrapping about the LBE core; and 
(b) the internal iteration slaved by the outer solver.

The present paper is the first part of this development
to focus on incorporating arbitrary nonlinear solution methods with the LBE. 
We begin with the primitive and weighted fixed-point iterations 
then carry on with
full approximation scheme (FAS) multigrid and Jacobian-free Newton--Krylov (JFNK) methods 
based on the same framework of the residual-driven LBE proposed in the present work.
The preconditioned LBE (PLBE) is further incorporated as an alternative nonlinear iteration and as a smoother within the FAS to demonstrate the ease of
solver composability. For the FAS, particular attention is paid to the coarse-grid defect construction and transfer operators compatible with the cell-centered nature of the LBE. A Chapman--Enskog analysis 
is used to show
the hydrodynamic consistency of the residual formulation and to identify the principal error components relevant to nonlinear convergence.

\begin{figure}[H]
    \centering
    \includegraphics[width=0.65\linewidth]{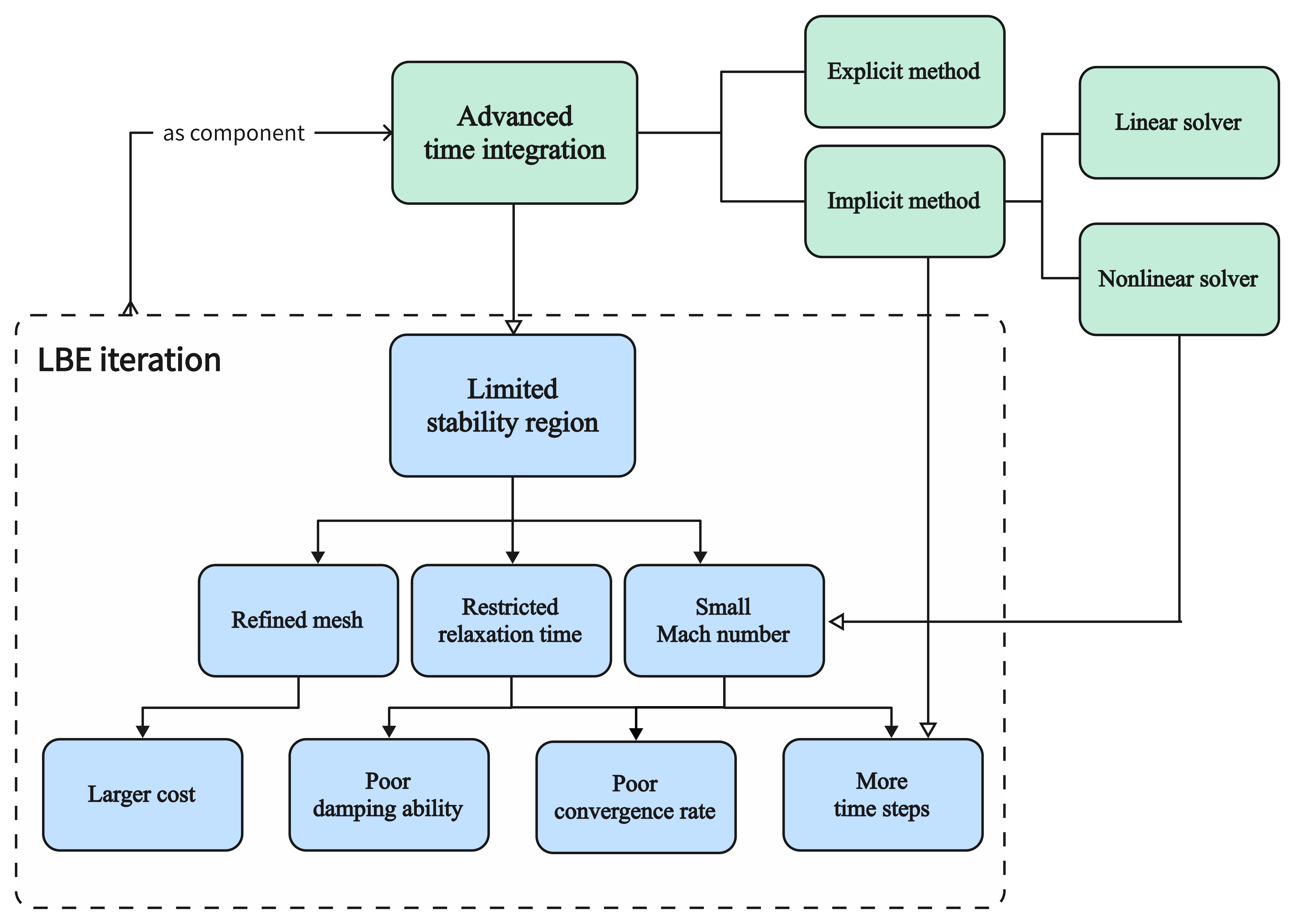}
    \caption{A schematic roadmap to wrap an LBE algorithm with an external iterative solver.}
\label{fig:roadmap}
\end{figure}

\cref{fig:roadmap} depicts a conceptual roadmap of our proposed residual-driven LBE 
wrapped with a nonspecific outer iteration solver.
The present work focuses on
the nonlinear-solver component, which provides the steady and pseudo-steady 
solution capability required by subsequent implicit physical-time formulations.
The success of this endeavor will pave the path for integrating other sophisticated 
time-stepping techniques with the LBE.

The remainder of this paper is organized as follows. 
\cref{sec:residual} introduces the residual-driven formulation of the LBE. The conventional collision-streaming update is expressed as a nonlinear discrete mapping, from which the distribution-function residual is constructed, and its relation to steady-state and physical-time formulations is discussed. 
\cref{sec:solvers} discusses a number of the nonlinear solution methods,
including fixed-point iteration, FAS multigrid, JFNK, and the combination with PLBE. 
\cref{sec:numerical} presents the results of numerical tests. 
The numerical tests, all in 2D, were chosen to achieve three objectives: 
(1) to verify that the nonlinear solvers do converge to the same discrete steady state of 
the primitive LBE; (2) to quantify the reduction in computational work achieved by the nonlinear solution strategies, and (3) to assess the composability of different nonlinear acceleration mechanisms through the common residual interface.
The numerical tests include (a) the lid-driven cavity flow at $\Rey=1000$ and 5000, (b) the 
flow past a cylinder at $\Rey=20$, and (c) the flow past a NACA0012 airfoil at 
an angle of attack $\alpha=5^\circ$ and $\Rey=1000$.
The cavity flow was used to verify fixed-point equivalence and quantify FAS efficiency; 
the flow past a cylinder was for assessing the accuracy with curved boundaries and 
the composition of FAS with PLBE; and the flow past a NACA0012 airfoil assesses the solver's 
capability for external flow in the low-Mach regime.
Finally, \cref{sec:con} summarizes the main findings and discusses
possible future extensions.
An Appendix discusses the consistency of the proposed residual-driven LBE method and 
error analysis.

\section{Residual-driven LBE}\label{sec:residual}

The central idea of the residual-driven LBE formulation is to delineate
the discrete LBE operator from a particular solution algorithm used for 
time advancing or otherwise.
One LBE collision-streaming step is encapsulated as a complete discrete update mapping, 
then a distribution-function residual is obtained directly from this mapping. 
The resulting residual obtained from the full set of distribution functions
provides the common input
for either steady-state computation or physical-time integration, while 
the entire LBE core algorithmic structure remains untouched and intact.

\subsection{Collision and Streaming process}

The LBE evolves
the distribution functions through a
sequence of local collision and streaming operations. 
The LBE algorithm core will be explicitly and succinctly 
discussed in what follows.
The key concepts and their associated symbols will be defined.
Considering the LBE on a \emph{cell-centered} uniform Cartesian grid 
in 2D with the grid spacing $\delta_x$:
\begin{equation}
    \{ \bm{r}_{\bm{j}} := 
	\bm{j} \delta_x,\ \bm{j} \in \mathbb{Z}_d  \},
    \quad d=2,
\end{equation}
and the velocity distribution functions on the grid:

\begin{equation}
    \qvector{f}
    =
    (f_0,\, f_1,\, \cdots,\, f_{q-1} )^\dagger,
    \quad
    f_\alpha :  \{ \bm{r}_{\bm{j}} \} \rightarrow \mathbb{R} .
\label{eqn:vectorf}
\end{equation}
where $\dagger$ indicates the transpose, so $\qvector{f}$ is a $q$-tuple column vector.
One LBE evolution step must consist of at least two essential steps --- collision and streaming:
\begin{subequations}
\label{eq:lbe-collision-streaming}
\begin{align}
    \text{Collision:\quad}
    &
    f^*_\alpha=f_\alpha
    + \delta_t \Omega_\alpha,
\label{eq:collision}\\
    \text{Streaming:\quad}
    &
    f_\alpha(\bm{r}_{\bm{j}}+\bm{\xi}_\alpha\delta_t,\, t_{m+1})
    =
    f^*_\alpha(\bm{r}_{\bm{j}},\, t_m).
\end{align}

\end{subequations}

The collision updates $f_\alpha$ \emph{locally}
by the 
process prescribed by $\delta_t\Omega_\alpha$, 
whereas streaming transports the
post-collision state $f^*_\alpha$ from one lattice node $\bm{r}_{\bm{j}}$
to a neighboring one $\bm{r}_{\bm{j}} + \delta_t \bm{\xi}_\alpha$,
determined by the corresponding particle velocity $\bm{\xi}_\alpha$. The
discrete particle velocities $\bm{\xi}_\alpha := \bm{c}_\alpha c$ constitute
the discrete velocity set $\mathbb{X}_q := \{ \bm{\xi}_\alpha | 0 \leq \alpha \leq (q-1) \}$
and $\mathbb{C}_q := \{ \bm{c}_\alpha \}$ is the lattice velocity set,
$c=\delta_x/\delta_t$ is the velocity unit.

The collision-streaming procedure described in \cref{eq:lbe-collision-streaming}, 
including the appropriate boundary conditions and effects due to forcing,
defines a nonlinear discrete mapping:
\begin{equation}
    \qvector{f}^{m+1}
    =
    \Tdelta (\qvector{f}^m),
\label{eq:T-map}
\end{equation}
where the superscript $m$ denotes successive time instances $t_m := m \delta_t$.
This process is conveniently interpreted as the evolution in one lattice
time step $\delta_t$.
In the residual-driven formulation, however, $\Tdelta$ is
regarded first as a discrete nonlinear mapping, while the physical-time
may be introduced separately and independently later
in \cref{sec:unified-residual}. 
That is, $\delta_t$ is the lattice time scale embedded in $\Tdelta$.

For the sake of concreteness, we will use the 2D LBE with nine discrete velocities, 
\emph{i.e.}, the D2Q9 model, of which the dimensionless lattice velocity set 
$\mathbb{C}_q$ includes:
\begin{align*}
\begin{array}{lllll}
        \bm{c}_0 = \begin{pmatrix}0\\0\end{pmatrix},
        &
        \bm{c}_1 = \begin{pmatrix}1\\0\end{pmatrix}, &
        \bm{c}_2 = \begin{pmatrix}0\\1\end{pmatrix}, &
        \bm{c}_5 = \begin{pmatrix}1\\1\end{pmatrix}, &
        \bm{c}_6 = \begin{pmatrix}-1\\1\end{pmatrix},
        \\
        &
        \bm{c}_3 =-\bm{c}_1,
        &
        \bm{c}_4 =-\bm{c}_2,
        &
        \bm{c}_7 =-\bm{c}_5,
        &
        \bm{c}_8 =-\bm{c}_6.
\end{array}
\end{align*}
In what follows, 
the BGK collision model will be used for the sake of simplicity:
\begin{equation}
  \delta_t \Omega_\alpha = - \omega [ f_\alpha-\fazero ],
\label{eqn:Omega}
\end{equation}
where $\omega$ is the dimensionless relaxation frequency, 
and the equilibrium distribution function $\fazero$ is
\begin{equation}
  \fazero(\hat{p},\, \bm{u})
  =
  \frac{w_\alpha}{\cs^2}
  \left[
    \hat{p} + \bm{\xi}_\alpha \cdot \bm{u}
    +\frac{1}{2}
    \left(
      \frac{(\bm{\xi}_\alpha \cdot \bm{u})^2}{\cs^2}
      -\bm{u}\cdot\bm{u}
    \right)
  \right] ,
\label{eqn:fzero}
\end{equation}
and the weight coefficients $w_\alpha$ for the D2Q9 LBE are
$$
	w_\alpha
	=
	\left\{
	\begin{array}{ll}
	4/9 , \quad & \alpha=0,
	\\
	1/9 , & 1 \leq \alpha \leq 4,
	\\
	1/36 , & 5 \leq \alpha \leq 8 .
	\end{array}
	\right.
$$
The kinematic viscosity of the LBE with the lattice BGK collision model
is:
\begin{equation}
  \nu = \cs^2 \delta_t 
  \left( \frac{1}{\omega}-\frac{1}{2}\right),
  \qquad
  \cs := \frac{c}{\sqrt{3}} ,
\end{equation}
where $\cs$ is the speed of sound in the LBE.

Although this work employs the lattice BGK (LBGK) model for the sake of simplicity,
it should be stressed that overall characteristics and qualitative features of
the proposed residual-driven LBE framework are independent of the specifics 
of the collision model. However, the specifics of the collision model may
affect numerical stability or convergence behavior. Since the focus of the present work is on establishing a general methodology for integrating nonlinear solvers with the LBE, the investigation of other collision models, such as multiple-relaxation-time (MRT) \cite{dHumieres2002ptrsla}, central-moment \cite{Geier2006pre}, and cumulant \cite{Geier2015cma} models, is left for future work.

\subsection{Distribution-function residual}
\label{sec:matrix-free-residual}

Mathematically, a complete collision-streaming process may be considered
as a particular means for nonlinear fixed point iteration.
Instead of restricting every application of
$\Tdelta$ to a physical-time advancement, we define its
defect as a nonlinear residual:
\begin{equation}
  \Adelta(\qvector{f})
  :=
  \frac{\qvector{f} - \Tdelta(\qvector{f})}{\delta_t},
\label{eq:A-operator}
\end{equation}
and the corresponding unscaled residual is
\begin{equation}
  \tensorfont{R}(\qvector f)
  :=
  \delta_t \Adelta(\qvector{f})
  =
  \qvector{f} - \Tdelta(\qvector{f}).
\label{eq:unscaled-residual}
\end{equation}
Both residuals,  $\Adelta(\qvector{f})$ and $\tensorfont{R}(\qvector{f})$,
are evaluated using the identical LBE evolution defined by $\Tdelta$.
Therefore,
the residuals themselves have no effect whatsoever on 
the underlying LBE evolution and its 
discretizations; they only provide an operational means
suitable for general nonlinear solution methods to be introduced.

A steady solution of the primitive LBE evolution is a fixed point of the
evolution process defined by $\Tdelta$, 
or equivalently it is the solution that nullifies 
the residual $\Adelta(\qvector{f})$:
\begin{equation}
\qvector{f}_* = \Tdelta(\qvector{f}_*)
  \Longleftrightarrow
  \Adelta(\qvector{f}_*) = \bm{0} .
\end{equation}
Thus, the primitive LBE iteration may be used as a
particular fixed-point iterative method for solving the nonlinear residual equation.
As will be discussed later, the FAS and Newton--Krylov methods also solve
the same equation, but only with
different residual-reduction strategies.

\begin{figure}[htb]
\centering
\includegraphics[width=0.3\linewidth]{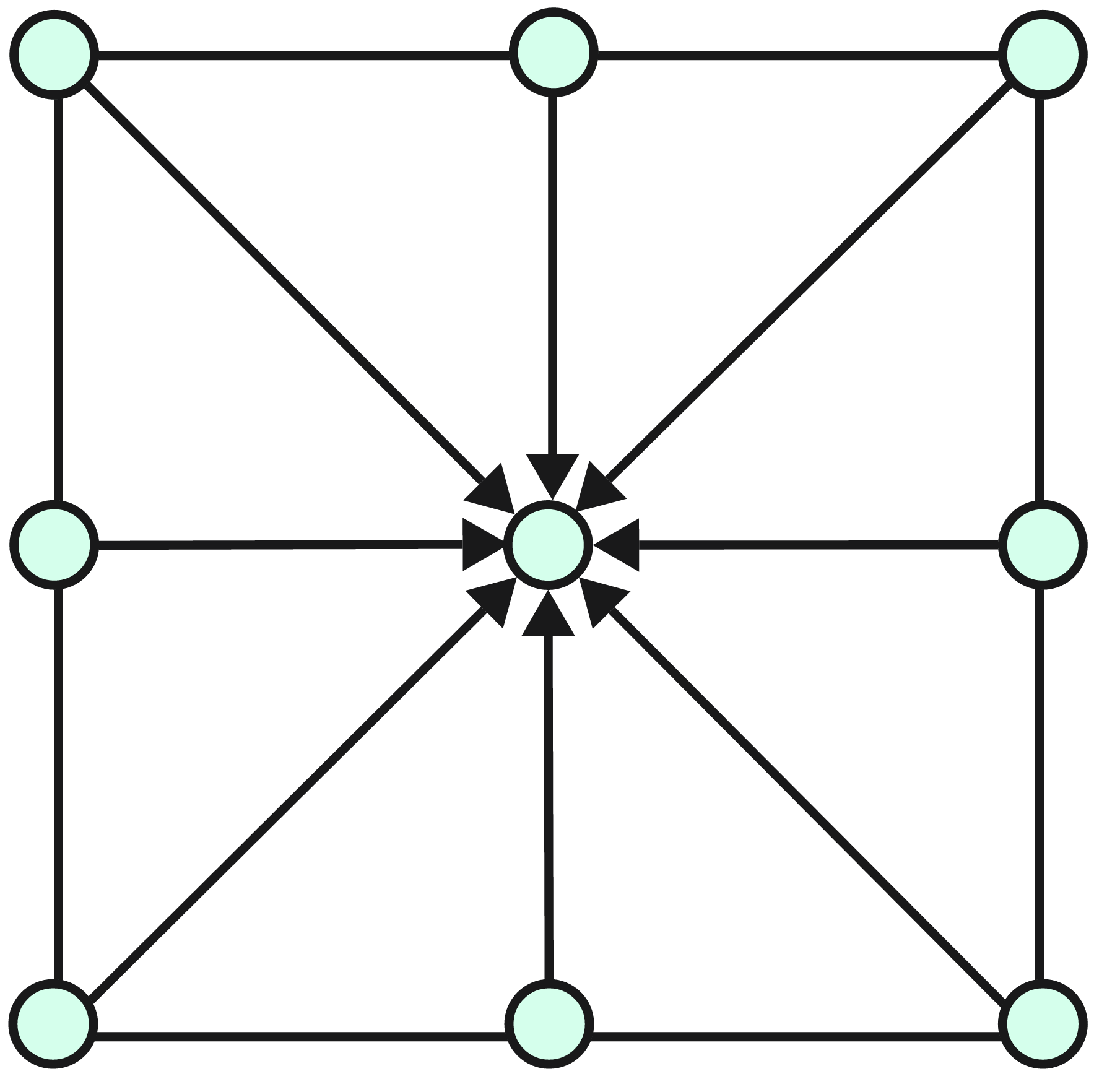}
\caption{Illustration of ``pull'' representation of the exact-streaming stencil.}
\label{fig:lattice}
\end{figure}

\subsection{Physical-time discretization}
\label{sec:unified-residual}

The residual equation \eqref{eq:unscaled-residual} can be used
for steady-state computation; in the meantime, the process of reaching a steady 
state is also driven by the LBE evolution process.
In other words, the process to reach a steady state is indeed ``physical'', 
that is, $\delta_t$ in the LBE system is a ``physical'' time step embedded in the 
evolution process. However, it is not necessary to be so --- a physical 
time step, $\delta_T$, may be introduced with an outer time integrator, 
independent of $\delta_t$ which is used within $\Tdelta$.
In this setting, $\Tdelta$ is acting as
an (inner) iterative means to achieve a steady state numerically, whether a 
temporary or perpetual one, so $\delta_t$ has no physical significance anymore and 
only becomes a pseudo time for pure numerical purposes. This demonstrates 
the possibility to separate the physical time $\delta_T$ dictated by an
outer time-integration scheme from the pseudo-time $\delta_t$ evolved with the
inner iteration, which is the LBE evolution process here. It should be
noted that advanced time-integration methods are not the focus of 
the present work; this brief discussion only serves as an outline. Instead,
general nonlinear solvers as the outer iterative techniques are to be studied 
here.

The residual operator $\Adelta$
may be used to introduce the following \emph{physical}-time evolution equation:
\begin{equation}
  \frac{\diff \qvector{f}}{\diff t}
  +\Adelta(\qvector{f})
  = \bm{0} ,
\label{eq:physical-evolution}
\end{equation}
where $t$ is the physical time now.  The spatio-temporal pair $(\delta_x,\, \delta_t)$ associated with
the LBE dynamics 
is embedded in the operator $\Adelta$: it determines exact node-to-node streaming and, together with the collision parameter $\omega$, determines the viscosity $\nu$.  Therefore $\delta_t$ is
not the step size selected by the \emph{outer} physical-time integrator to treat the time-evolution of $\diff\qvector{f}/ \diff t$.
The particular time integrators introduced below are used only to illustrate the 
underlying mechanism in the proposed framework to separate the physical-time step $\delta_T$
from the lattice time step $\delta_t$ embedded in the LBE operator.

\paragraph{Explicit physical-time integration}

\cref{eq:physical-evolution} can be advanced explicitly with
a standard explicit ODE method. Such a method evaluates $\Adelta$ with only a given known state 
$\qvector{f}^m$
or at a specific
time-level $m$
such that it requires no nonlinear inner iteration.  
For example, a second-order explicit Runge-Kutta (RK) 
method may be utilized as the following:
\begin{subequations}
\label{eq:rk2-explicit-lbe}
\begin{align}
  \qvector{f}^{\star}
  &=\qvector{f}^n - \delta_T \,
    \Adelta(\qvector f^n), \\
  \qvector f^{n+1}
  &=\qvector f^n-\frac{\delta_T}{2}
    \left[
      \Adelta(\qvector f^n)
      +\Adelta(\qvector f^{\star})
    \right] ,
\end{align}
\end{subequations}
alternatively, the two-step Adams-Bashforth scheme is:
\begin{equation}
  \qvector{f}^{n+1}
  =\qvector{f}^n - \frac{\delta_T}{2}
   \left[
     3 \Adelta(\qvector{f}^n)
     - \Adelta(\qvector{f}^{n-1})
   \right].
\label{eq:ab2-explicit-lbe}
\end{equation}
In the aforementioned explicit schemes based on the residual-driven LBE framework,
$\delta_T$ is the physical time step, while the lattice time step
$\delta_t$ is hidden
within the LBE dynamics, only relevant 
to the residual evaluation. In this setting,
$\delta_T$ need not be connected
with $\delta_t$ physically or otherwise, although its admissible size 
is determined and restricted by the stability region of the specific method
to discretize the time-derivative $\diff \qvector{f}/\diff t$ as well as the spectrum of the discrete LBE residual operator.
In this setting,
the primitive LBE is fully recovered as the forward-Euler scheme, the
simplest member in this family of explicit schemes, when $\delta_T = \delta_t$:
\begin{equation}
  \qvector{f}^{n+1}
  =
  \qvector{f}^n - \delta_t \Adelta(\qvector{f}^n)
  = 
  \Tdelta(\qvector{f}^n) .
\label{eq:forward-euler-primitive-lbe}
\end{equation}
The above consideration broadens the LBE's numerical horizons concerning its time-integration.

\paragraph{Implicit physical-time integration}
The explicit schemes outlined above demonstrate that separating the physical time step size
$\delta_T$ from the lattice time step $\delta_t$ immediately introduces 
the freedom to choose physical-time discretization without any hindrance of the spatial
resolution $\delta_x$ defined by the LBE spatial discretization. 
Nevertheless, explicit integration schemes still bear the birthmark of the
primitive LBE algorithm, because its 
stability is determined in part by the spectrum of the discrete LBE residual operator, 
which could still be very restrictive. 
But, much larger physical time steps
may be allowed when an implicit scheme is employed as the outer iteration with the 
proposed residual-driven LBE framework.

Denote $\delta_T := T_{n+1}-T_n$ the \emph{physical} time step.
A fully implicit discretization of \cref{eq:physical-evolution} requires the solution of
\begin{equation}
  \tensorfont{G}^{n+1}(\qvector{f})
  :=
  \tensorfont{D}_{T}(\qvector{f})
   + \Adelta(\qvector{f})
  = \bm{0} ,
\label{eq:general-residual}
\end{equation}
for the unknown state $\qvector{f} := \qvector{f}^{n+1}$.  
For an $M$-step formula, the temporal term may be written as
\begin{equation}
  \tensorfont{D}_{T}(\qvector{f})
  = \frac{a_0 \qvector{f}-\qvector{E}(\qvector{f}^n,\, \ldots,\,
    \qvector{f}^{n+1-M})}{\delta_T},
\label{eq:general-time-operator}
\end{equation}
where $\qvector{E}$ is a known linear combination of previous $M$ physical-time states.
The BDF2 scheme is a particular example:
\begin{equation}
  \tensorfont{G}^{n+1}(\qvector{f})
  = 
  \frac{3 \qvector{f}-4 \qvector{f}^n+\qvector{f}^{n-1}}{2\delta_T}
   +\Adelta(\qvector{f})
  =\bm{0} ,
\label{eq:bdf2-residual}
\end{equation}
where the nonlinear residual equation is solved implicitly, 
rather than with a sequence of primitive LBE steps over the interval $\delta_t$.

A pseudo-time variable $\tau$ can be introduced solely for the purpose 
to drive the nonlinear residual to zero:
\begin{equation}
  \frac{\partial\qvector{f}}{\partial\tau}
  +\tensorfont{G}^{n+1}(\qvector{f})
  =
  \bm{0} .
\label{eq:pseudo-time-equation}
\end{equation}
The variable $\tau$ has no physical significance;
an increment in $\tau$ only indicates an advance in the numerical solution process and
does not advance the physical solution from $T_n$ to $T_{n+1}$.  
When converged, \eqref{eq:pseudo-time-equation} leads to
\begin{equation}
  \frac{\partial \qvector{f}}{\partial\tau}
  \longrightarrow
  \bm{0} ,
  \qquad
  \tensorfont{G}^{n+1}(\qvector{f}) \longrightarrow \bm{0} ,
\end{equation}
so that the converged pseudo-steady state is the implicit solution 
$\qvector{f}^{n+1}$ at the physical time $T_{n+1}$.

In a recent work \cite{JiLuo2026DTSLBE}, 
the local physical-time term at the new pseudo iteration $(k+1)$ 
and the residual at the current iteration $k$ are combined as the following:
\begin{equation}
  \frac{\qvector{f}^{k+1}-\qvector{f}^{k}}{\delta_\tau}
  +\tensorfont{D}_{T}(\qvector{f}^{k+1})
  +\Adelta(\qvector{f}^{k})
  =
  \bm{0} ,
\label{eq:pseudo-time-split}
\end{equation}
where $\delta_\tau$ is a solver parameter, \emph{i.e.}, the pseudo-time step.
For the BDF2 scheme of \eqref{eq:bdf2-residual},
\cref{eq:pseudo-time-split} becomes:
\begin{equation}
  \left(1+\frac{3\delta_\tau}{2 \delta_T}\right) \qvector{f}^{k+1}
  =
  \qvector{f}^{k}
   -\delta_\tau \Adelta (\qvector{f}^{k})
   + \frac{\delta_\tau}{2 \delta_T}
      \left(4\qvector{f}^n-\qvector{f}^{n-1} \right).
\label{eq:bdf2-pseudo-iteration}
\end{equation}

For a steady state, the physical-time term is absent,
\begin{equation}
  \tensorfont{D}_{T} = \qvector 0,
  \qquad
  \tensorfont{G}(\qvector f) = \Adelta(\qvector{f})
  =
  \bm{0} .
\label{eq:steady-residual}
\end{equation}
Then \cref{eq:pseudo-time-split} with $\delta_\tau = \delta_t$ reduces to the primitive LBE fixed-point iteration.  Again, $\delta_t$ plays a role in both the LBE dynamics and in the
relaxation process of the steady-state solver, but not as the physical time step 
in the evolution process. In contrast, the explicit updates in 
\cref{eq:rk2-explicit-lbe,eq:ab2-explicit-lbe,eq:forward-euler-primitive-lbe} 
do advance the physical solution in time directly and involve no pseudo-time variable.
Therefore they ought not be confused with the pseudo-time iterations used to solve
an implicit equation.

\subsection{Existing implicit LBE schemes}
\label{sec:relation-dts}
There have been implicit LBE schemes \cite{tolke1998implicit,huang2015fully,li2019conservation},
which do not necessarily construct their nonlinear residual equations from the same 
native LBE discrete operators. This essential distinction is
not merely about the choice of a particular nonlinear solver,
but about the residual to which the nonlinear solver is applied.
\cref{tab:implicit-lbe-conceptual} compiles some existing implicit LBE 
approaches according to this distinction.

\begin{table*}[htbp]
\centering
\caption{Features of some existing implicit kinetic schemes related to the LBE
}
\label{tab:implicit-lbe-conceptual}
\footnotesize
\setlength{\tabcolsep}{3pt}
\renewcommand{\arraystretch}{1.18}

\begin{tabular}{
p{0.25\textwidth}
p{0.65\textwidth}}
\toprule
\textbf{Formulation}
&
\textbf{Discrete operator}
\\
\midrule

Implicit FD
\cite{tolke1998implicit,huang2015fully}
and FV \cite{li2019conservation} kinetic schemes
&
The DVBE is discretized directly with either FD 
or FV methods, leading to implicit schemes based on DVM. 
The nonlinear solver acts on the residual of the
FD/FV discretization of the DVBE.
\\

\addlinespace

MG dual-time-stepping LBE \cite{gsell2020multigrid}
&
The primitive LBE algorithm is used as the pseudo-time iteration. 
The physical-time correction is defined with only the conserved-moments.
\\

\addlinespace

Residual-driven implicit LBE in \cite{JiLuo2026DTSLBE}
&
The complete collision-streaming map
$\Tdelta$ is retained as the discrete LBE operator. The implicit equation is built from the residual of the standard LBE itself.
\\

\bottomrule
\end{tabular}
\end{table*}

The finite-difference \cite{tolke1998implicit,huang2015fully}
and finite-volume \cite{li2019conservation}
kinetic methods discretize the discrete-velocity Boltzmann equation (DVBE)
\cref{eqn:dvbe0}, so the derived 
discretizations differ from those native to the LBE \cite{Junk2001nmpde}.
Other than using the same set of discrete velocities, 
their connections to the LBE are therefore tenuous otherwise; thus they are
the implicit DVM schemes (cf., \emph{e.g.}, \cite{dimarco2027cicp}).
We noted the LBE \emph{per se} is 
a very specific one-step, second-order splitting discretization of the
DVBE \cite{Dellar2013cma}, and its spatial discretizations are indeed somewhat
unconventional \cite{Junk2001nmpde}. The favorable features of the LBE
discretizations
as well as the efficiency associated with 
its simplicity have been lost in these implicit schemes
\cite{tolke1998implicit,huang2015fully,li2019conservation}, while they
are fully retained in the methodology proposed in this work.

The existing multigrid dual-time-stepping method \cite{gsell2020multigrid} 
does retain the basic collision-streaming algorithmic core as its pseudo-time 
iteration, 
however its 
physical-time residual is derived only
from the conserved moments, as opposed to the full set of the distribution functions, 
and is included into the kinetic equations as source terms.
In contrast, the residual-driven
formulation considered here instead defines
$$
    \Adelta(\qvector{f})
    =
    \frac{\qvector{f}-\Tdelta(\qvector{f})}{\delta_t}
$$
directly from the complete LBE evolution process.
Consequently, both steady problems and implicit physical-time stages are formulated
in terms of the residual of the LBE itself.

The dual-time-stepping LBE in a recent work
\cite{JiLuo2026DTSLBE} is a BDF2 realization of this
distribution-level residual formulation, which employed
MG-LBE as the nonlinear inner solver and demonstrated the feasibility
of implicit physical-time integration. The present study does not
focus on time-integration methods. Instead, we investigate the
efficiency, robustness, and composition of general nonlinear solvers for the
steady and implicit residual equations 
related to the LBE.

\subsection{Hydrodynamic consistency and convergence criteria}
\label{sec:consistency-tolerance}
The preceding discussions are about 
the residual equations to be solved and 
their relations to existing implicit LBE schemes.
Before introducing nonlinear solution methods, 
it remains to determine 
the accuracy required by the solutions of the residual equations.

The solution of \cref{eq:general-residual} driven by the LBE with the 
low-Mach-number scaling should converge to the leading-order hydrodynamic
equations dictated by the primitive LBE in the limit of $\delta_t \rightarrow 0$; 
the solution process of \cref{eq:general-residual} certainly induces 
errors caused by physical-time discretization and nonlinear iteration, 
in addition to those in the primitive LBE algorithm. Therefore, the total error
may be written as the following:
\begin{equation}
  E_{\mathrm{total}}
  =O(\delta_t^2)
   +O(\delta_T^\kappa)
   +O(\zeta_{\ssstyle \mathrm{NL}})
   +O(\delta_t \May^2),
\label{eq:error-decomposition}
\end{equation}
where $\kappa$ is the order of the outer physical-time method, and
$\zeta_{\ssstyle \mathrm{NL}}$ denotes the solution error due to an inexact
nonlinear solve. This solution error is distinct from the residual tolerance
$\varepsilon_{\ssstyle \mathrm{NL}}$ defined below. 
A bound on
$\zeta_{\ssstyle \mathrm{NL}} \leq C_{\ssstyle \mathrm{NL}}\varepsilon_{\ssstyle \mathrm{NL}}$
requires a suitable local residual-to-solution stability estimate; the
constant includes the normalization used in the halting criterion. 
The lattice and physical-time truncation terms are examined through
Chapman--Enskog analysis in the Appendix, while controlling the solution
error due to an inexact nonlinear solver additionally requires a
residual-to-solution stability estimate.
\cref{eq:error-decomposition} will be used
to identify the error components relevant to 
inner-solver tolerances in what follows.

The nonlinear iteration is terminated when the normalized residual
reaches its given tolerance:
\begin{equation}
  \frac{\|\tensorfont{G}^h(\qvector{f}^{(k)})\|_2}
       {\max\{\|\tensorfont{G}^h(\qvector{f}^{(0)})\|_2,1\}}
  \leq \varepsilon_{\ssstyle \mathrm{NL}}.
\label{eq:nonlinear-stopping}
\end{equation}
For a steady problem, this criterion monitors $\Adelta$. For an implicit solution, 
the complete residual $\tensorfont{G}$, including the temporal term, must be monitored. The criterion alone does not ensure a solution-error bound without residual-to-solution stability.

The nonlinear equations with a common convergence criterion are now specified, 
the numerical methods used to reduce these residuals efficiently will be discussed
next.

\section{Nonlinear solvers}
\label{sec:solvers}

With the residual equation and its convergence criterion given explicitly,
different nonlinear algorithms can now be compared fairly. We begin with the weighted fixed-point iteration, which provides both a natural baseline and the nonlinear smoother used within multigrid. We then develop an FAS multigrid scheme and a Newton--Krylov scheme for the same residual, and finally consider compositions of these solvers. 
In what follows,
the discrete LBE operator itself remains unchanged; only the strategy to reduce the
residual toward zero is changed.

\subsection{Fixed-point iteration}
\label{sec:fixed-point}

The weighted residual iteration,
\begin{equation}
  \qvector{f}^{k+1}
  =\qvector{f}^k-\lambda\delta_t
  \Adelta(\qvector{f}^k) ,
\label{eq:weighted-fixed-point}
\end{equation}
reduces to a primitive LBE update with $\lambda=1$.
This iterative scheme is used as the baseline and as the nonlinear smoother on each FAS level. Its low cost is advantageous, but its convergence deteriorates with mesh refinement and may 
slow down significantly 
when $\omega \rightarrow 2^-$.

As an alternative single-grid iteration, we also consider the
preconditioned lattice Boltzmann equation (PLBE)
\cite{guo2004preconditioned}. For the D2Q9 model, the corresponding
equilibrium distribution is:
\begin{equation}
\fazero
=
\frac{w_\alpha}{c_s^2}
\left[
\hat p
+
\bm{\xi}_\alpha \cdot \bm{u}
+
\frac{1}{2\gamma}
\left(
\frac{(\bm{\xi}_\alpha\cdot\bm{u})^2}{c_s^2}
-
\bm{u}\cdot\bm{u}
\right)
\right].
\end{equation}
The standard LBE is recovered for $\gamma=1$. For steady states, the
PLBE and standard LBE recover the same macroscopic
steady equations, whereas their transient eigenvalue systems differ.
The preconditioning parameter
$\gamma \leq 1$ modifies the equilibrium distribution so as to reduce
the disparity between the characteristic propagation speeds of the
underlying kinetic system.
Consequently, choosing $\gamma<1$ can substantially accelerate
convergence to a steady state without modifying the targeted 
steady-state solution. Therefore, the PLBE can be used
as an alternative nonlinear
iteration for the same steady residual equation rather than a different
target problem. The relaxation time $\tau_p$ of the PLBE is related to that
of the corresponding standard LBE by
\begin{equation}
    \frac{1}{\omega}-\frac{1}{2}
    =
    \gamma\left(\frac{1}{\omega_p}-\frac{1}{2}\right).
\end{equation}

\subsection{FAS multigrid}
\label{sec:fas}
 
The primitive LBE fixed-point iteration becomes increasingly inefficient with mesh refinement due to its slow damping of long-wavelength hydrodynamic error modes.
To overcome this limitation, we employ the full approximation scheme (FAS) which acts directly on the nonlinear residual rather than on a linearized correction equation. Its construction for the present LBE requires three ingredients: (1) a consistent coarse-grid defect equation, 
(2) nonlinear smoothing and coarse-grid operators, 
and (3) transfer operators compatible with the cell-centered 
discretizations and its boundary treatment.
Implementation details of the FAS may be found in existing work 
\cite{mavriplis2006multigrid,gsell2020multigrid,JiLuo2026DTSLBE}.

\subsubsection{FAS defect for the LBE}

For the steady unscaled residual $\tensorfont{R}^h=\delta_t^h \tensorfont{A}^h$,
the defect is
\begin{equation}
  \qvector{D}^H
  =\tensorfont{R}^H(\tensorgreek{\Phi}_p\qvector{f}^h)
   -\frac{\delta_t^H}{\delta_t^h}
    \tensorgreek{\Phi}_b \tensorfont{R}^h(\qvector{f}^h),
\label{eq:implemented-defect}
\end{equation}
where $\tensorgreek{\Phi}_p$ and $\tensorgreek{\Phi}_b$ are the state and
the residual restriction operators, respectively; their definitions
will be given in \cref{subsec:transfer} later. 
The superscripts $H$ and $h$ indicate coarse and fine grids, respectively.
The factor $\delta_t^H/\delta_t^h$ follows from the use of unscaled residuals, 
so it is not an empirical correction.

\subsubsection{Nonlinear smoothing and coarse-grid operators}

The steady nonlinear smoother is
\begin{equation}
  \qvector{f}^{h,k+1}
  =\qvector{f}^{h,k}
   -\lambda
   \left[\tensorfont{R}^h(\qvector{f}^{h,k})-\qvector{D}^h\right].
\label{eq:fas-smoother}
\end{equation}
The same relaxation frequency $\omega$ is used on all levels. 
Under factor-two coarsening with $\delta_t^H=2\delta_t^h$ and an identical
lattice speed $c = \delta_x/\delta_t$ (which amounts to the acoustic scaling),
the viscosities in two grids are related as the following:
\begin{equation}
  \nu^H=2\nu^h
\end{equation}
for a fixed relaxation frequency $\omega$. The coarse equation is therefore an auxiliary 
defect equation rather than a physical simulation with the same Reynolds number.

The hierarchy that preserves physical viscosity at all grid levels
would instead satisfy
\begin{equation}
  \frac{1}{\omega^H}-\frac{1}{2}
  =\frac{\delta_t^h}{\delta_t^H}
   \left(\frac{1}{\omega^h}-\frac{1}{2}\right).
\label{eq:constant-viscosity-omega}
\end{equation}
For factor-two coarsening,
$\omega^H=4\omega^h/(2+\omega^h)$, which compels
$\omega^H$ closer to $2$ as the grid is coarsened and
may weaken the smoothing property. The fixed-$\omega$ hierarchy is therefore 
merely an algorithmic choice.

This restriction on $\omega^H$
can be relaxed when PLBE is used as the nonlinear
smoother. Owing to the following relationship:
$$
    \frac{1}{\omega_p}-\frac{1}{2}
    =
    \frac{1}{\gamma} 
	\left( \frac{1}{\omega}-\frac{1}{2} \right) ,
$$
the PLBE relaxation frequency $\omega_p$ remains away from the limiting value
$\omega_p=2$ for $\gamma<1$. This allows the relaxation
frequency $\omega$ to be rescaled between successive grid levels while
preserving the target physical Reynolds number $\Rey$ as a constant. 
Accordingly, for FAS-PLBE we use
\begin{equation}
    \frac{1}{\omega_p^H}-\frac{1}{2}
    =
    \frac{\delta_t^h}{\delta_t^H}
    \left(\frac{1}{\omega_p^h}-\frac{1}{2}\right),
\end{equation}
so that the same Reynolds number is maintained at all grid levels.

\subsubsection{Cell-centered transfer and boundary treatment}
\label{subsec:transfer}

In the present LBE discretization, including various boundary conditions,
the distribution functions reside on
cell centers. The multigrid transfer operators must therefore be defined on a cell-centered hierarchy, rather than 
the vertex-centered operators commonly used in geometric 
multigrid \cite{mohr2004cell,JiLuo2026DTSLBE}. As a consequence, the fine- and coarse-grid nodes are 
no longer collocated as in the vertex centered grids. 
We use piecewise-constant restriction for the state and bilinear restriction for the residual:
\begin{equation}
 \tensorgreek{\Phi}_p=\frac{1}{4}
  \begin{bmatrix}
    1&1\\
    1&1
  \end{bmatrix},
  \qquad
 \tensorgreek{\Phi}_b=\frac{1}{64}
  \begin{bmatrix}
    1&3&3&1\\
    3&9&9&3\\
    3&9&9&3\\
    1&3&3&1
  \end{bmatrix}.
  \label{eq:transfer-stencils}
\end{equation}
Following the cell-centered transfer construction of
Mohr and Wienands \cite{mohr2004cell}, the coarse-grid correction is
prolonged using the
operator $\tensorfont{P}_b$ defined as the adjoint
of the residual-restriction operator $\tensorgreek{\Phi}_b$. 
The values outside of a Dirichlet boundary are set to zero when bilinear transfer is evaluated.

\subsubsection{Multigrid cycle}

\begin{figure}[H]
\centering
\begin{minipage}[b]{0.30\textwidth}
\centering
\text{(a) V-cycle}

\vspace{0.3em}

\begin{tikzpicture}[scale=0.58, every node/.style={font=\scriptsize}]
    \draw[gray!50] (0,0) -- (3.8,0);
    \draw[gray!50] (0,-1) -- (3.8,-1);
    \draw[gray!50] (0,-2) -- (3.8,-2);

    \node[left] at (0,0) {$L_1$};
    \node[left] at (0,-1) {$L_2$};
    \node[left] at (0,-2) {$L_3$};

    \draw[thick]
        (0,0) -- (0.95,-1) -- (1.9,-2)
        -- (2.85,-1) -- (3.8,0);
\end{tikzpicture}
\end{minipage}
\hfill
\begin{minipage}[b]{0.34\textwidth}
\centering
\text{(b) W-cycle}

\vspace{0.3em}

\begin{tikzpicture}[scale=0.58, every node/.style={font=\scriptsize}]
    \draw[gray!50] (0,0) -- (4.8,0);
    \draw[gray!50] (0,-1) -- (4.8,-1);
    \draw[gray!50] (0,-2) -- (4.8,-2);

    \node[left] at (0,0) {$L_1$};
    \node[left] at (0,-1) {$L_2$};
    \node[left] at (0,-2) {$L_3$};

    \draw[thick]
        (0,0) -- (0.8,-1) -- (1.6,-2)
        -- (2.4,-1) -- (3.2,-2)
        -- (4,-1) -- (4.8,0);
\end{tikzpicture}
\end{minipage}
\hfill
\begin{minipage}[b]{0.30\textwidth}
\centering
\text{(c) F-cycle}

\vspace{0.3em}

\begin{tikzpicture}[scale=0.58, every node/.style={font=\scriptsize}]
    \draw[gray!50] (0,0) -- (4,0);
    \draw[gray!50] (0,-1) -- (4,-1);
    \draw[gray!50] (0,-2) -- (4,-2);

    \node[left] at (0,0) {$L_1$};
    \node[left] at (0,-1) {$L_2$};
    \node[left] at (0,-2) {$L_3$};

    \draw[thick]
        (0,-2) -- (1,-1) -- (2,-2)
        -- (3,-1) -- (4,0);
\end{tikzpicture}
\end{minipage}
\caption{Schematic representation of the multigrid cycle types.
$L_1$, $L_2$, and $L_3$ denote successively coarser grid levels.}
\label{fig:mg_cycles}
\end{figure}
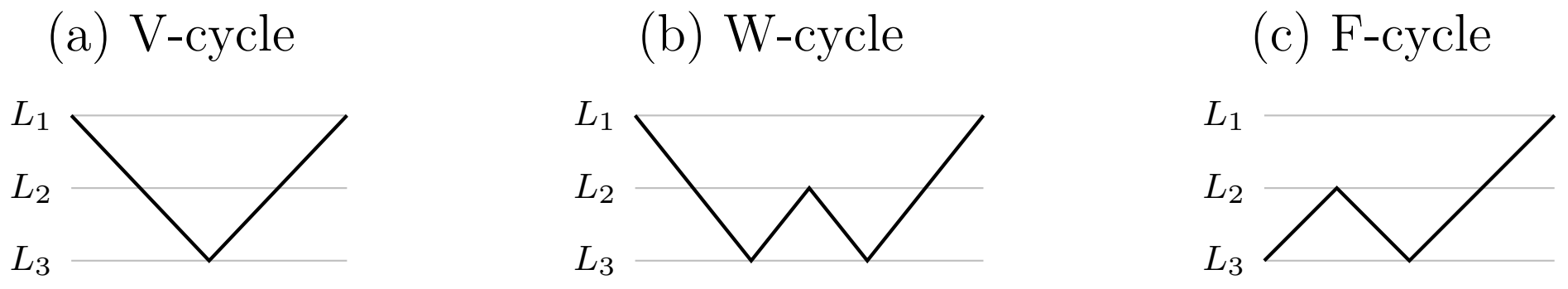

The preceding components can now be assembled into a complete nonlinear multigrid iteration.
\cref{fig:mg_cycles} illustrates
the V-, W-, and F-cycle traversal patterns used in the multigrid solver. The V- and W-cycles correspond to $\mu=1$ and $\mu=2$ in Algorithm~1, respectively, whereas the F-cycle employs the coarse-to-fine traversal shown in the figure. \cref{alg:fas} summarizes the resulting V- and W-cycles in terms of pre-smoothing, coarse-grid defect correction, prolongation, and post-smoothing.

\begin{algorithm}[H]
\caption{FAS-LBE V-cycle ($\mu=1$) or W-cycle ($\mu=2$)}
\label{alg:fas}
\begin{algorithmic}[1]
\Function{FAS}{$\qvector{f}_\ell,\qvector D_\ell,\ell$}    \Comment{$\ell$ denotes the \textit{l}-th level}
  \If{$\ell$ is the coarsest level}
    \State \Return $\mathrm{Smooth}^{\nu_3}(\qvector{f}_\ell,\qvector{D}_\ell)$
  \Else
    \State $\qvector{f}_\ell^{\mathrm{pre}}\gets
      \mathrm{Smooth}^{\nu_1}(\qvector{f}_\ell,\qvector{D}_\ell)$
    \State $\qvector{f}_{\ell-1}^{0}\gets \tensorgreek{\Phi}_p\qvector{f}_\ell^{\mathrm{pre}}$
    \State Form $\qvector{D}_{\ell-1}$ using \cref{eq:implemented-defect}
    \State $\qvector{f}_{\ell-1}\gets\qvector{f}_{\ell-1}^{0}$
    \For{$i=1,\ldots,\mu$}
      \State $\qvector{f}_{\ell-1}\gets\Call{FAS}{\qvector{f}_{\ell-1},\qvector{D}_{\ell-1},\ell-1}$
    \EndFor
    \State $\widetilde{\qvector{f}}_\ell\gets\qvector{f}_\ell^{\mathrm{pre}}
      +\tensorfont{P}_b(\qvector{f}_{\ell-1}-\qvector{f}_{\ell-1}^{0})$
    \State \Return $\mathrm{Smooth}^{\nu_2}(\widetilde{\qvector{f}}_\ell,\qvector{D}_\ell)$
  \EndIf
\EndFunction
\end{algorithmic}
\end{algorithm}

\subsection{Newton--Krylov method}
\label{sec:newton}

The FAS achieves acceleration by transferring nonlinear error across grid scales. 
The Newton--Krylov methods treat the same problem differently
by constructing a local linear model for the nonlinear residual
then solving for a correction. The residual-driven formulation makes this comparison particularly transparent, since the Newton Jacobian is exactly
the derivative of 
the same LBE
residual used by both the fixed-point iteration and FAS.

For the steady nonlinear equation
\begin{equation}
  \tensorfont{R}(\qvector{f})=
  \bm{0} ,
\end{equation}
the Jacobian of the LBE residual is
\begin{equation}
  \tensorfont{J_R}(\qvector{f})
  :=
  \frac{\partial \tensorfont{R}}{\partial \qvector{f}}(\qvector{f})
  =
  \tensorfont{I} -
  \frac{\partial \Tdelta}{\partial \qvector{f}}(\qvector{f}) .
  \label{eq:residual-jacobian}
\end{equation}
A damped Newton step is obtained from
\begin{equation}
  \tensorfont{J_R}(\qvector{f}^k)\,\delta\qvector{f}^k
  =-\tensorfont{R}(\qvector{f}^k),
  \qquad
  \qvector{f}^{k+1}=\qvector{f}^k+\eta_k\delta\qvector{f}^k,
\label{eq:exact-newton}
\end{equation}
where $\eta_k$ is selected by damping or line search. The linear system in \cref{eq:exact-newton} may be solved inexactly by a Krylov method.
Importantly, the Newton operator corresponds to the derivative of the actual
LBE residual $\tensorfont{R}$, rather than
a surrogate one. The efficiency depends on the quality of the linear solver, preconditioning, and globalization, while its local convergence behavior differs substantially from that of FAS.

In the Jacobian-free implementation, the Jacobian matrix is not assembled
explicitly. Instead, its action on a Krylov vector $\qvector{v}$ is approximated
by a directional finite difference:
\begin{equation}
    \tensorfont{J_R}(\qvector{f})\qvector{v}
    \approx
    \frac{
    \tensorfont{R}(\qvector{f}+\varepsilon\qvector{v})-\tensorfont{R}(\qvector{f})
    }{\varepsilon},
\end{equation}
where $\varepsilon$ is a finite-difference perturbation parameter.
Consequently, the Krylov iteration requires only evaluations of the same
LBE residual used by the fixed-point and FAS solvers,
without forming or storing the full Jacobian matrix
\cite{KnollKeyes2004}. An approximate Jacobian may nevertheless be used
as a preconditioner for the Krylov iteration.

A Jacobian-free Newton--Krylov (JFNK) method \cite{KnollKeyes2004} is implemented 
here, but its detailed optimization is not a primary focus of the present work. It is
introduced mainly to demonstrate that the residual-driven formulation
provides a common interface through which different nonlinear solvers can
be applied and combined without modifying the underlying LBE
discretization.

\section{Numerical Results}
\label{sec:numerical}

The numerical results are organized progressively from a canonical
internal-flow problem to flows with curved boundaries, and then external-flow configurations.
The 2D lid-driven cavity is used first to verify the discrete fixed point and to
investigate the convergence characteristics and efficiency of FAS. The flow past
a cylinder at $\Rey=20$ examines the accuracy of the method for a flow problem 
involving a curved boundary and nontrivial inflow/outflow boundary conditions, while also
assessing the effectiveness and composability of FAS with preconditioning. 
Finally, the flow past a NACA0012 airfoil is considered to
demonstrate the applicability of the proposed framework to steady external aerodynamic
flow and to examine its efficiency in the low-Mach regime. All nonlinear
solvers are implemented within PETSc \cite{petsc-user-ref,Brune2015}.
All wall-clock comparisons were performed under the same computational
environment and are reported only as relative implementation-level references.

\subsection{2D lid-driven cavity flow}

\paragraph{Fixed-point equivalence and solution accuracy}
We first verify numerically that the nonlinear solution procedure does not
alter the discrete steady problem.  The cavity with a mesh size of $128^2$ 
at $\Rey=1000$ is solved independently with the primitive LBE iteration, FAS-LBE, and JFNK-LBE. The fully converged primitive LBE solution,
denoted by $\qvector{f}_{\ssstyle \mathrm{LBE}}$, is used as the reference.
\cref{tab:fixed_point} compiles the relative global $L_2$ differences in $\qvector{f}$
between the nonlinear-solver solutions and the LBE reference:
\begin{equation}
	\varepsilon
	:=
	\frac{\sum_j \|\qvector{f}(\bm{r}_j) - \qvector{f}_{\ssstyle \mathrm{LBE}} (\bm{r}_j)\|_2}
         {\sum_j \|\qvector{f}_{\ssstyle \mathrm{LBE}}(\bm{r}_j)\|_2} .
\end{equation}

\begin{table}[H]
\centering
\caption{2D lid-driven cavity flow at $\Rey=1000$ with $128^2$ mesh.
  Comparison of fixed-point solutions with the primitive LBE reference.}
\label{tab:fixed_point}
\begin{tabular}{lr}
\toprule
    \multicolumn{1}{c}{Method}
    &
	\multicolumn{1}{c}{$\varepsilon$}
    \\
    \midrule
    FAS-LBE      & $5.5\times10^{-13}$ \\
    JFNK-LBE     & $3.8\times10^{-13}$ \\
    \bottomrule
  \end{tabular}
\end{table}

The differences with the two nonlinear solution methods are close, and they
both are at the level of $10^{-13}$, \emph{i.e.}, the level of numerical 
round-off. This result attests the fixed-point equivalence of the residual 
formulation --- the FAS and Newton methods only alter the nonlinear solution
procedure without affecting the underlying LBE being solved.

\begin{figure}[!htbp]
\centering
\includegraphics[width=0.80\linewidth]{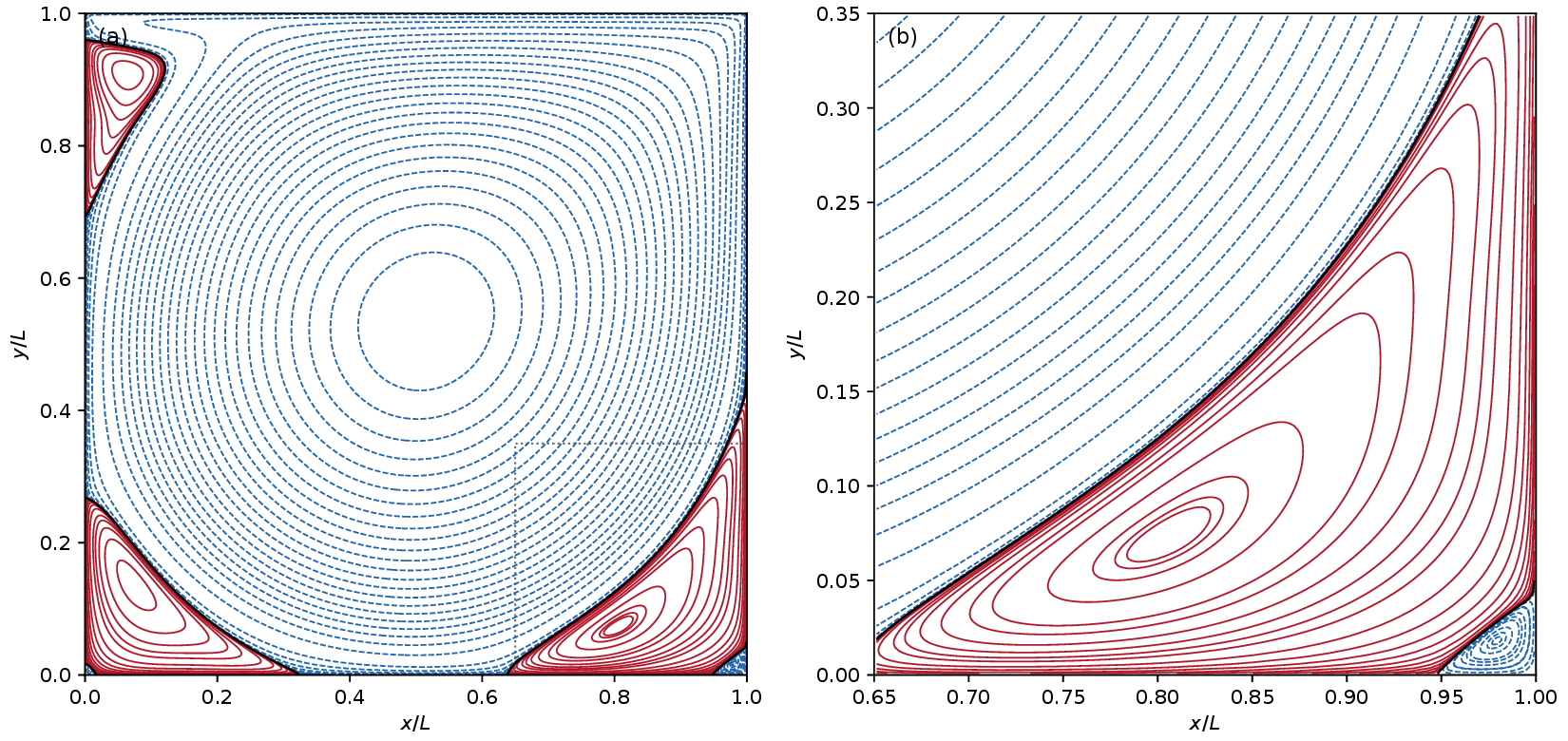}
\caption{2D lid-driven cavity flow at $\Rey=5000$ with mesh size
  $512^2$. The streamlines obtained with the FAS-LBE scheme.
}
\label{fig:re5000_stream}
\end{figure}

The flow at $\Rey=5000$ with a larger mesh size of $512^2$, which is a much
more challenging case, is used to assess the accuracy of the converged flow 
fields. \cref{fig:re5000_stream} shows the streamline obtained with the FAS-LBE.  
The figure clearly shows that the secondary vortex near the upper left corner, and
a tertiary corner vortex in the lower right corner, have all been captured.
The tertiary vortex is particularly difficult for the primitive LBE scheme 
to capture \cite{HouSL1995jcp}.

\begin{figure}[!htbp]
\centering
\includegraphics[width=0.9\linewidth]{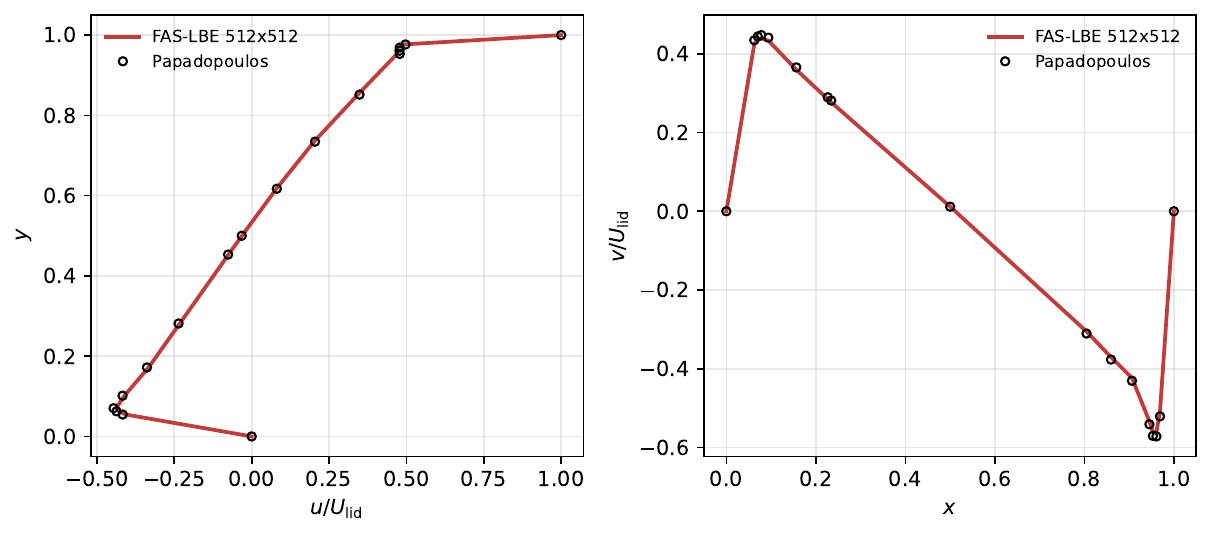}
\caption{2D lid-driven cavity flow at $\Rey = 5000$ with mesh size of $512^2$.
The normalized centerline velocities 
$u(x=L/2,\, y)/U_0$ (left) and $v(x,\, y=L/2)/U_0$ (right)
obtained with the FAS-LBE scheme compared to the data
in \cite{Papadopoulos2015}.}
\label{fig:centerline_velocity}
\end{figure}

\begin{figure}[htbp]
\centering
\includegraphics[width=0.90\linewidth]{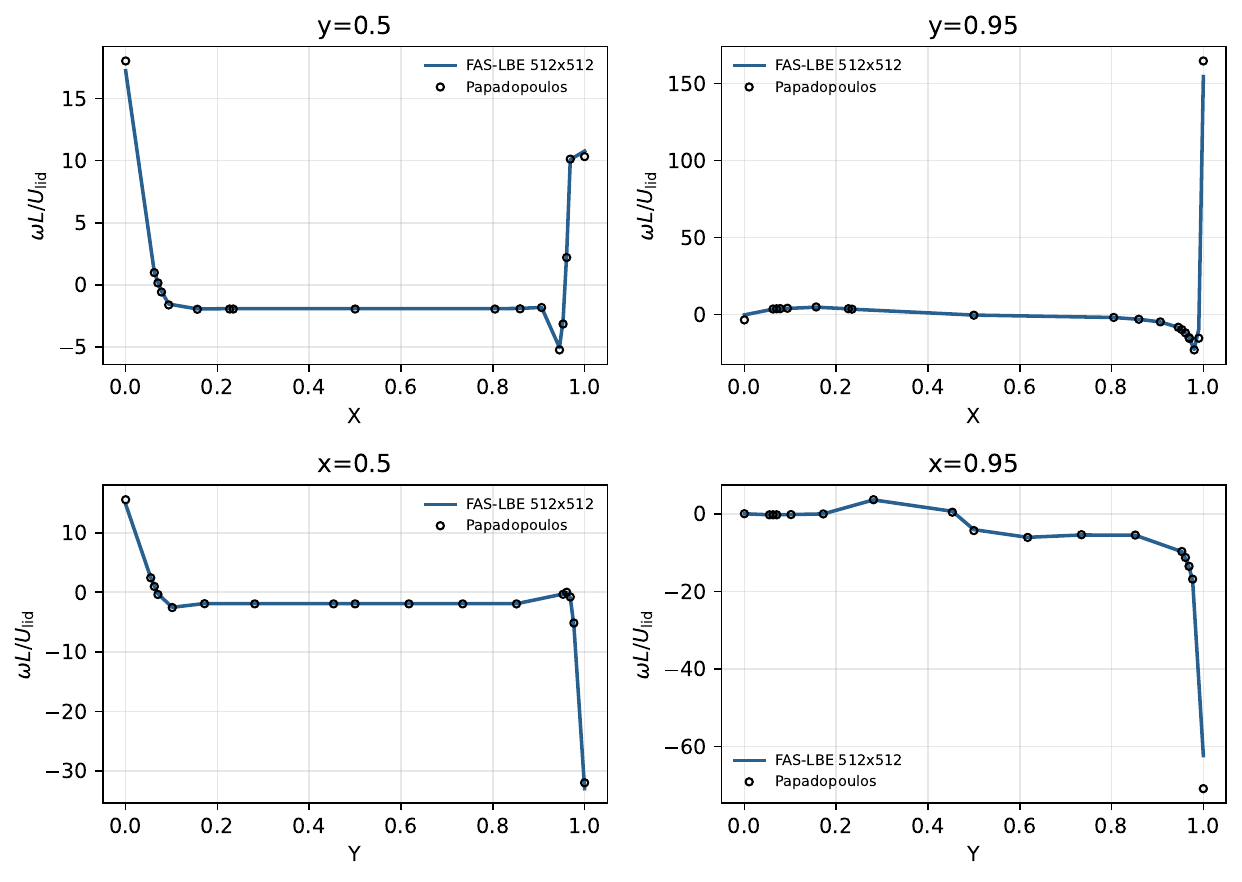}
\caption{2D lid-driven cavity flow at $\Rey=5000$, mesh size of $512^2$. The normalized 
vorticity $\varpi(x, y) L/U_0$ along two horizontal lines,
$y/H=0.5$ and 0.95 (top row), and two vertical ones $x/L=0.5$ and 0.95 (bottom row),  
compared with the data in \cite{Papadopoulos2015}.}
\label{fig:vorticity_validation}
\end{figure}

Further quantitative comparisons of the FAS-LBE solution at $\Rey=5000$ 
are made in \cref{fig:centerline_velocity} and \cref{fig:vorticity_validation}
with the solutions in \cite{Papadopoulos2015}, which were obtained
with the second-order FD time-marching solution based on the vorticity-stream 
function formulation and Richardson extrapolations to reach converged 
solutions. The meshes used in  \cite{Papadopoulos2015} were as large as 
$1661^2$. Our results of the normalized centerline velocities
$u(x=L/2,\, y)/U_0$ and $v(x,\, y=L/2)/U_0$, where $U_0$ is the lid velocity, differed from those in \cite{Papadopoulos2015} by 1.8\%;
in the figures, they are nearly indistinguishable.

The vorticity field, $\varpi := \bm{\nabla} \times \bm{u}$,
may provide a more sensitive measure of the near-wall solution quality.
\cref{fig:vorticity_validation} compares the normalized vorticity field $\varpi L/U_0$ at the centerlines, 
near the moving lid at $y/H=0.95$
and near the right wall at $x/L=0.95$. 
The FAS-LBE solutions agree very well with the high-quality
data in \cite{Papadopoulos2015}
overall, with the exceptions at the points immediately adjacent to the
moving wall near top right corner,
where the solution is singular and 
the vorticity field is particularly sensitive to the 
resolution and wall treatment. It should be noted that the mesh size used in 
\cite{Papadopoulos2015} was much larger than ours used here, and the data
in \cite{Papadopoulos2015} were obtained with Richardson extrapolations on
the data computed from three large meshes.
The comparisons do illustrate the fact that the FAS-LBE retains the accuracy of
the underlying LBE discretization.

\paragraph{Computational efficiency}
Once both the fixed-point equivalence and the accuracy of the
converged steady solution are verified, 
the computational work required
to reach accurate solutions is ready to be assessed.
The FAS configurations used in the two principal
timing comparisons are listed
in \cref{tab:fas_timing_parameters},
where $\nu_1$ and $\nu_2$ denote the numbers of pre- and post-smoothing
sweeps, respectively, and $\nu_3$ denotes the number of coarse-grid
smoothing sweeps (cf. Algorithm \ref{alg:fas}).

\begin{table}[htbp]
  \centering
  \caption{FAS parameters used for the principal FAS-LBE timing
  comparisons.}
  \label{tab:fas_timing_parameters}
  \begin{tabular}{cccccccc}
    \toprule
    Mesh size & $\Rey$ & $\lambda$ & $\nu_1$ & $\nu_2$ & $\nu_3$
    & cycle & levels \\
    \midrule
    $128^2$ & 1000 & 0.8 & 3 & 3 & 20 & W & 5 \\
    $512^2$ & 5000 & 0.6 & 5 & 5 & 40 & W & 7 \\
    \bottomrule
  \end{tabular}
\end{table}

To further quantify the computational work of a multigrid cycle independent
of wall-clock time, we define the work units (WU) in terms of
the cost in the units of one finest-grid LBE iteration,
\emph{i.e.}, one sweep on the finest-grid.
For an $L$-level W-cycle with $\nu_1$ pre-smoothing and $\nu_2$
post-smoothing sweeps on each non-coarsest level, and $\nu_3$ iterations
on the coarsest level, the theoretical WU per cycle is therefore
\begin{equation}
W^*
=
(\nu_1+\nu_2)
\sum_{\ell=0}^{L-2} 2^{-\ell}
+
\nu_3 2^{-(L-1)}=2(\nu_1+\nu_2)+ \frac{\nu_3-2(\nu_1+\nu_2)}{2^{L-1}}.
\label{eq:wcycle_work}
\end{equation}

\cref{tab:wcycle_work} compares the idealized estimates
given by \cref{eq:wcycle_work} with the
measured equivalent WU of one FAS W-cycle, $\overline{W}$. The measured values are
slightly larger because Eq.~\eqref{eq:wcycle_work} accounts only for the
LBE smoothing, but ignores other operations in an actual FAS cycle including
residual evaluations, restriction and prolongation, construction of the
coarse-grid FAS defect, and other multigrid operations.  Nevertheless,
the measured costs remain close to the theoretical estimates, indicating
that these auxiliary operations constitute only a moderate fraction of
the total cycle cost, about 20\%.

\begin{table}[htbp]
\centering
\caption{Theoretical and measured WU of one FAS W-cycle, expressed
in equivalent finest-grid LBE iterations.}
\label{tab:wcycle_work}
\begin{tabular}{cccccc}
\toprule
Levels & $(\nu_1,\nu_2,\nu_3)$
& $W^*$
& $\overline{W}$
& $\overline{W}/ W^*$
\\
\midrule
5 & $(3,\, 3,\, 20)$ & 12.5000 & 14.9000 & 1.192 \\
7 & $(5,\, 5,\, 40)$ & 20.3125 & 22.8373 & 1.124 \\
\bottomrule
\end{tabular}
\end{table}

\begin{figure}[!htbp]
\centering
\includegraphics[width=0.75\linewidth]{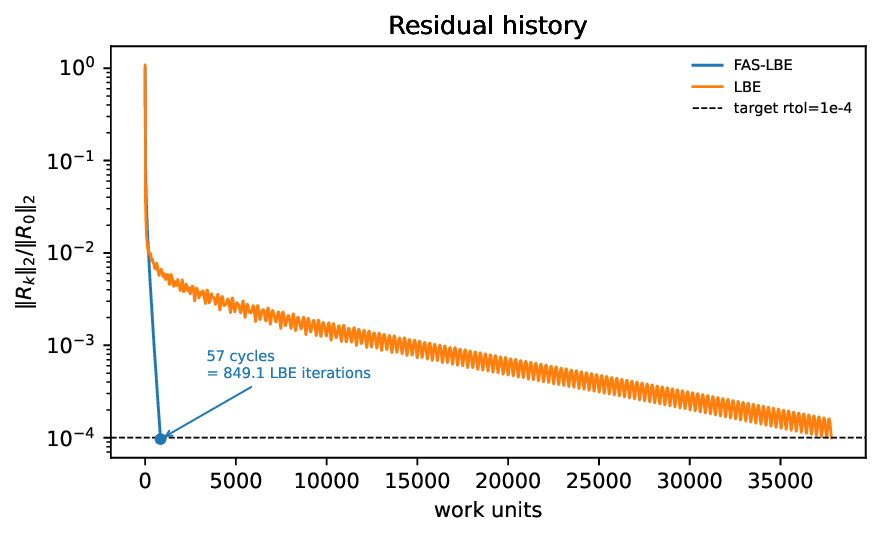}
\caption{2D lid-driven cavity flow at $\Rey=1000$ with mesh size $128^2$.
The evolution of the normalized residuals in the FAS-LBE and 
the primitive LBE iteration in terms of WU.
}
\label{fig:fas_lbe_re1000}
\end{figure}

\cref{fig:fas_lbe_re1000} compares the evolution of the residual
in the FAS-LBE with that in the primitive LBE iteration
for the 2D lid-driven cavity flow with $\Rey=1000$
and mesh size of $128^2$.
Both calculations halted after the nonlinear residuals were reduced by
four orders of magnitude, \emph{i.e.},
$\| \qvector{R}_k \|_2 / \| \qvector{R}_0 \|_2 \le 10^{-4}$.
The FAS-LBE reaches this tolerance in 57 cycles, equivalent to
approximately 849 LBE iterations on the finest mesh, 
and requires about $1.04$s of measured solver time.  The primitive LBE requires
$46.2$s.  The acceleration due to the FAS is therefore
about $46.2/1.04 \approx 44.4$ in terms of wall clock times.

\begin{figure}[!htbp]
\centering
\includegraphics[width=0.75\linewidth]{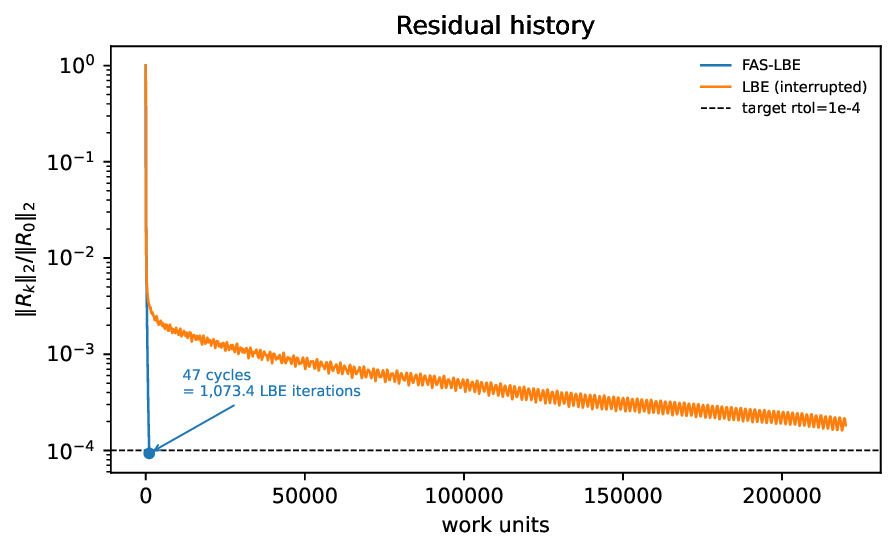}
\caption{2D lid-driven cavity flow at $\Rey=5000$ and mesh size of $512^2$.
The evolution of the residuals with the FAS-LBE and the primitive LBE 
in terms of WU.
}
\label{fig:fas_lbe_re5000}
\end{figure}

The efficiency gain or computational speedup becomes substantially greater for
the case with $\Rey=5000$ and a larger mesh size of $512^2$.
As \cref{fig:fas_lbe_re5000} shows, the FAS-LBE
reaches the same $10^{-4}$ relative-residual tolerance 
in 47 cycles, equivalent to
approximately 1073 LBE iterations on the finest mesh, and requires
$27.9$s. In contrast, the primitive LBE iteration halted
after $2.2 \cdot 10^5$ iterations, which took approximately $6250$s. 
The gain in computational speed due to the FAS is 224 in this case.
The comparisons above clearly demonstrated the advantages of the nonlinear multigrid 
solver in terms of efficiency. This efficiency gain is particularly pronounced when the primitive LBE direct fixed-point iteration becomes inefficient.

\paragraph{Comparison with Newton--Krylov and a hybrid solver}
We have repeatedly stressed the fact that the proposed residual-driven framework
commands the freedom to accommodate and incorporate \emph{any} nonlinear solver without 
affecting the inherent LBE algorithmic kernel. To demonstrate this point, a pure 
Jacobian-free Newton--Krylov (JFNK) solver and a simple hybrid scheme which applies one FAS cycle
before switching to JFNK are implemented to compare with the FAS solver.
A comparison is made with the cavity flow at $\Rey=1000$ and mesh size of 
$128^2$. The nonlinear tolerance to halt is $10^{-8}$.

The FAS scheme took
only 4.01s, whereas the JFNK scheme cost about
507.77s, and the hybrid FAS-JFNK scheme consumed about
379.93s. This unequivocally demonstrates the superiority
of the FAS scheme; its time is only about 
0.8\% of that of the JFNK scheme. In fact,
applying one FAS iteration before JFNK reduces the time by approximately 25\%.
This justifies the choice of the FAS scheme as the primary nonlinear solver in what follows.

A more meaningful comparison is made with the finite-element method (FEM)
using Taylor-Hood elements implemented in \textsc{FreeFEM} platform \cite{freefem}.
The FreeFEM uses $P_2$ and $P_1$
elements for the velocity and pressure field, respectively,
on unstructured meshes. Both the FAS-LBE and FEM 
were implemented with \textsc{PETSc} for the nonlinear/linear algebra infrastructure.
\cref{tab:freefemsetup} provides the details for the \textsc{FreeFEM} implementation.
The test case is the cavity flow at $\Rey=1000$. For the FAS-LBE,
the mesh contains $128^2$ lattice nodes, corresponding
to $9\times128^2=147,456$ distribution-function unknowns for D2Q9.
The FEM uses 32,258 triangular elements and has 146,433 independent
algebraic degrees of freedom (DOFs). Thus, the two implementations have
comparable numbers of algebraic unknowns, rather than comparable mesh-node counts.
The relative nonlinear residual tolerance is set at $10^{-4}$ for the test.
The times for the FAS-LBE vs. the FEM are 1.02s vs. 16.01s, so the FAS-LBE
is about 15 times faster than the FEM for this particular case.

\begin{table}[htbp]
\centering
\caption{FreeFEM/PETSc configuration}
\label{tab:freefemsetup}
\begin{tabular}{ll}
\toprule
Components
& Setting \\
\midrule
Mesh & 32,258 triangular elements \\
Independent DOFs & 146,433 \\
Nonlinear solver & PETSc Newton line search \\
Linear solver & PETSc FGMRES \\
Preconditioner & LU/MUMPS \\
\bottomrule
\end{tabular}
\end{table}

\paragraph{Dependence of convergence on grid resolution $N$, $\omega$, and $\lambda$}

Finally we investigate
the dependence of the FAS convergence on
the grid resolution $N$, the relaxation frequency $\omega$
(cf. \cref{eqn:Omega}), and the weight $\lambda$ 
(cf. Eqs. \eqref{eq:weighted-fixed-point} 
and \eqref{eq:fas-smoother}).
The multigrid scheme in the FAS is a W-cycle with five levels,
and $\nu_1=\nu_2 = 3$, $\nu_3=20$.

With $\lambda = 0.6$, we consider two sequences of grid-refinement with a given
grid Reynolds number $\Rey^* := \Rey/N$. The first one is:
$$
(N,\, \Rey ) = (128,\, 1000),\; (256,\, 2000),\; (512,\, 4000),
\; \Rey^* = 7.8125;
$$
and the second one is:
$$
(N,\, \Rey )=(128,\, 100),\;(256,\, 200),\;(512,\, 400) ,
\; \Rey^* = 0.78125 .
$$
A fixed $\Rey^*$ simply means that $\omega$ remains as a constant
when the resolution $N$ varies.

\begin{figure}[!htbp]
\centering
\begin{subfigure}[t]{0.475\textwidth}
\centering
\includegraphics[width=\linewidth]{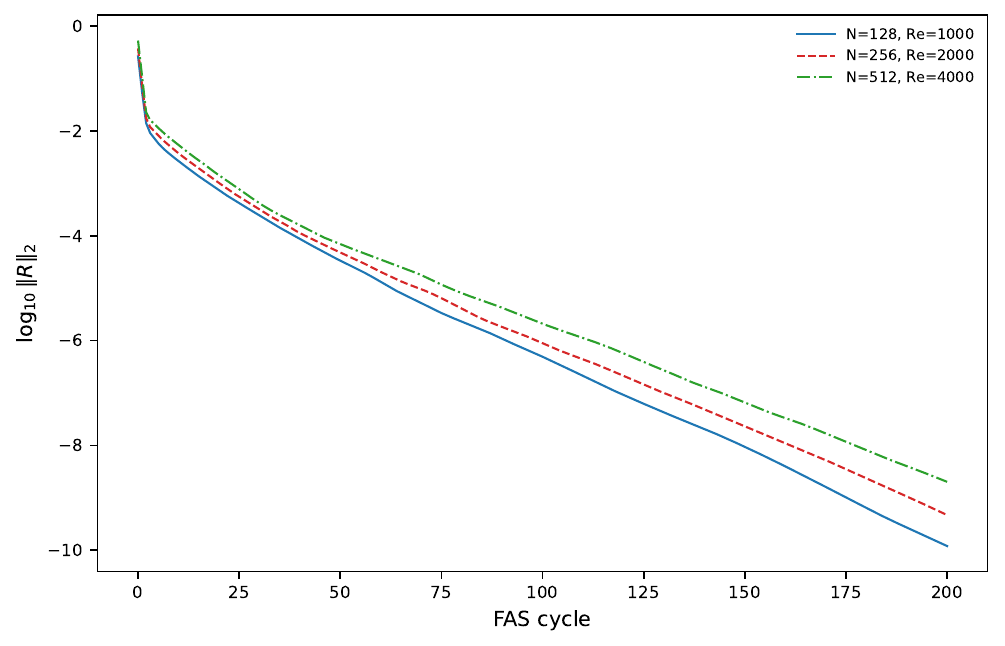}
\caption{$(N,\Rey)=(128,1000),(256,2000),(512,4000)$}
\label{fig:fas_high}
\end{subfigure}
\hfill
\begin{subfigure}[t]{0.475\textwidth}
\centering
\includegraphics[width=\linewidth]{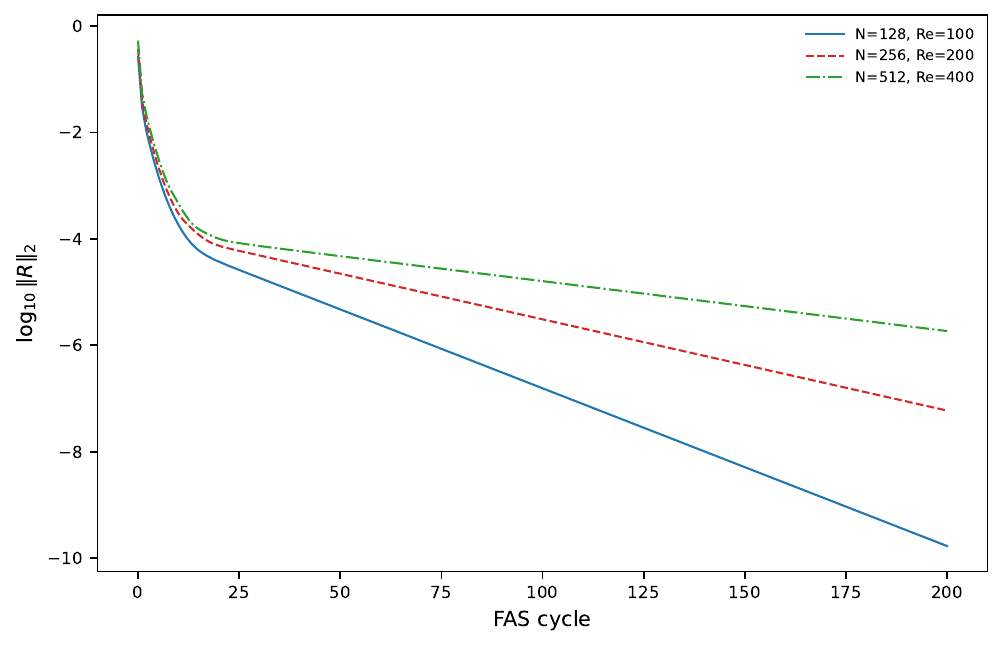}
\caption{$(N, \Rey)=(128,100),(256,200),(512,400)$}
\label{fig:fas_low}
\end{subfigure}

\vspace{0.5em}

\begin{subfigure}[htbp]{0.95\textwidth}
\centering
\includegraphics[width=\linewidth]{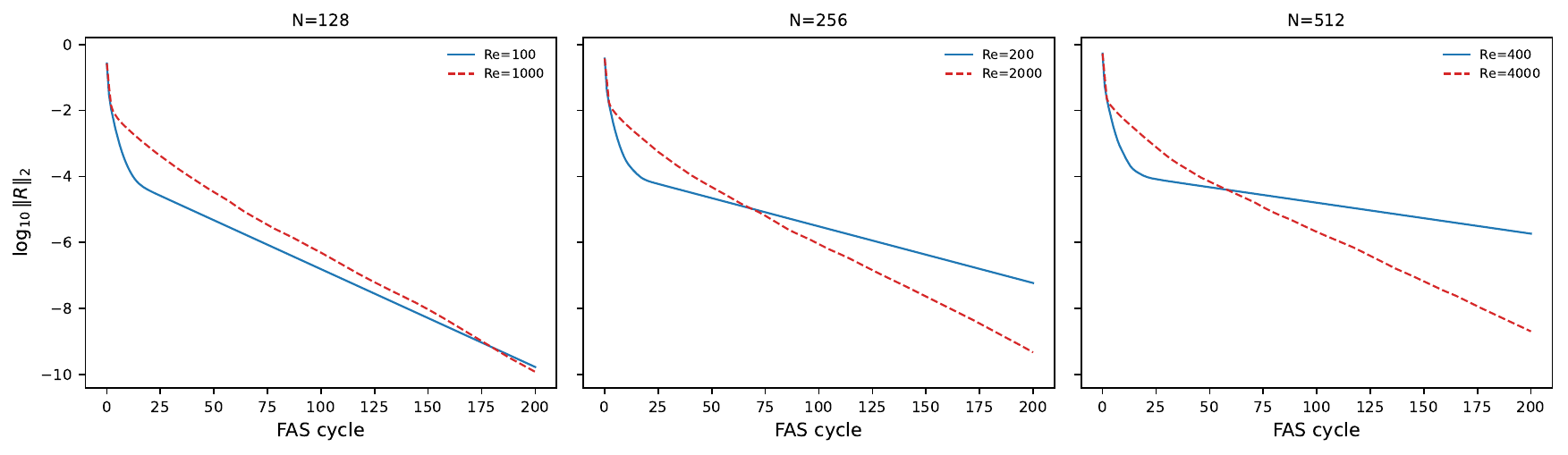}
\caption{Comparisons at fixed grid resolution $N$.}
\label{fig:fas_pairs}
\end{subfigure}
\caption{The evolution of the FAS-LBE residual with $\lambda=0.6$
}
\label{fig:fas_parameter_dependence}
\end{figure}

\cref{fig:fas_parameter_dependence} shows
the evolution of the residuals obtained with different parameter settings
aforementioned. The following observations are to be made. First of all, with the 
higher value of $\Rey^* = 7.8125$, the asymptotic convergence rate does 
not seem to depend on the resolution $N$ much; $N$ only has a very 
modest effect on the number of iterations to attain convergence,
as clearly shown in \cref{fig:fas_parameter_dependence}(a).
A markedly different behavior is observed for the lower value 
of $\Rey^* = 0.78125$: after the initial rapid reduction, the asymptotic 
convergence rate becomes progressively slower as the grid is refined, 
\emph{i.e.}, as $N$ increases, as shown in \cref{fig:fas_parameter_dependence}(b).

The comparison may be reorganized for a fixed $N$, and varying $\Rey^*$, as shown
in \cref{fig:fas_parameter_dependence}(c). The results represented in this figure
show that, with $N$ fixed, the asymptotic rate of convergence for the residual increases 
with the $\Rey^*$ or $\Rey$, indicating the effect of $\omega$. This effect is
particularly conspicuous for the case of $N=512$, in which the residual for
the higher $\Rey^*$ case decays much faster than the lower $\Rey^*$ counterpart.
This observation may be counter-intuitive to the experience that a steady state
at a high $\Rey^*$ may be more difficult to solve than its lower $\Rey^*$
counterparts, which is certainly true for the primitive LBE algorithm.
Therefore, we may conclude that while the FAS iteration has substantially 
mitigated the grid-dependent slowdown associated with the primitive LBE algorithm,
it cannot completely remove it. That is, the 
asymptotic convergence rate of the residual does depend on $\omega$,
which is intrinsic to the LBE kernel. This is not surprising because 
the LBE algorithm is used for the nonlinear smoothing and coarse-grid 
operations (cf. details in~\cref{sec:fas}).

\begin{figure}[!htbp]
    \centering

\begin{subfigure}[t]{0.475\textwidth}
\centering
\includegraphics[width=\linewidth]{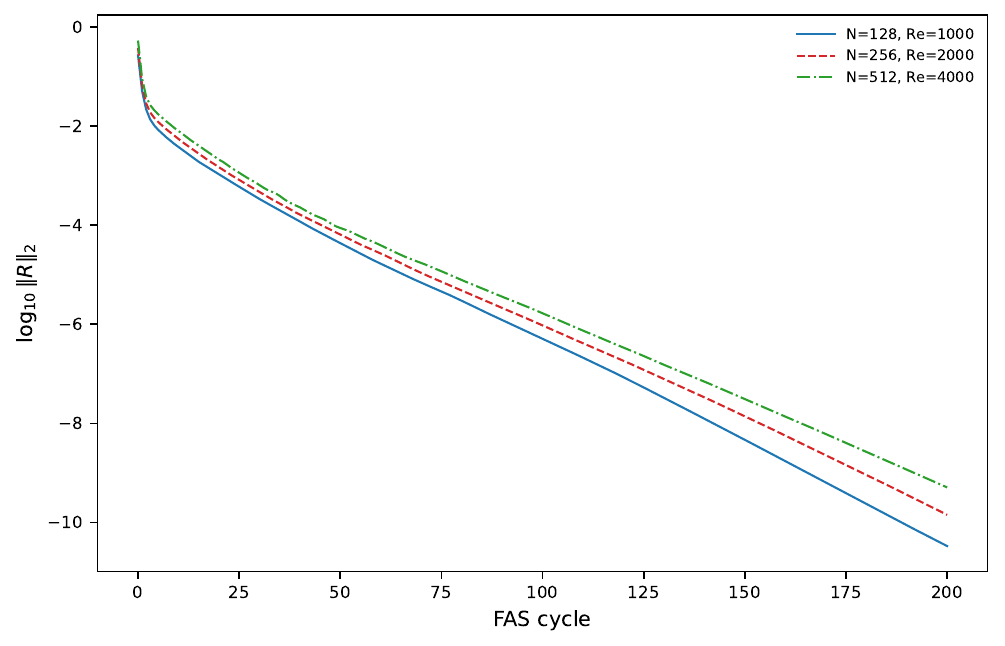}
\caption{$(N, \Rey)=(128,1000),(256,2000),(512,4000)$.}
\label{fig:fas_high_08}
\end{subfigure}
\hfill
\begin{subfigure}[t]{0.475\textwidth}
\centering
\includegraphics[width=\linewidth]{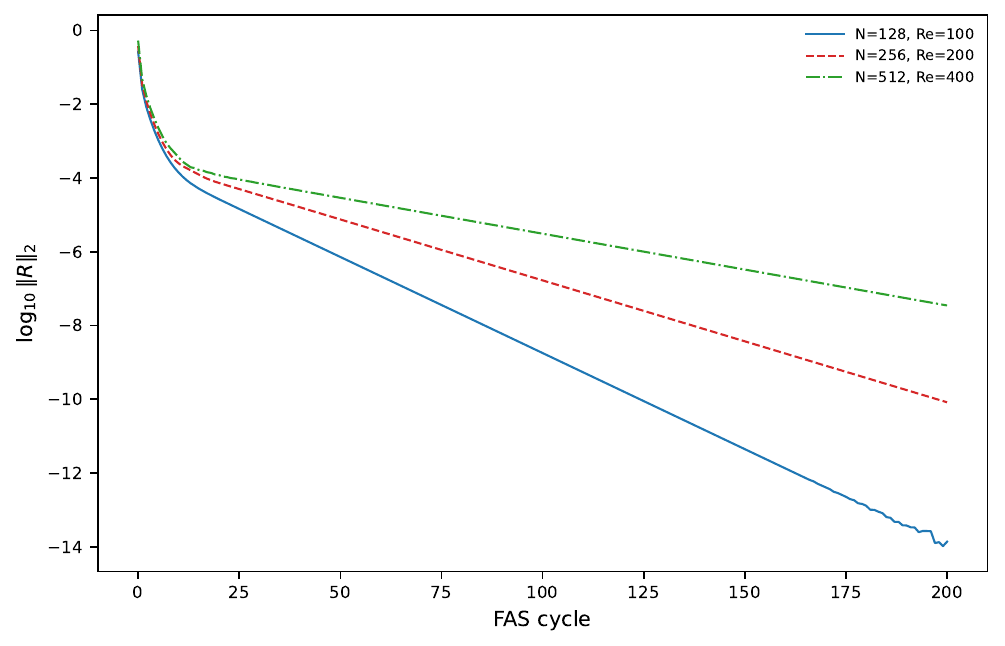}
\caption{$(N, \Rey)=(128,100),(256,200),(512,400)$.}
\label{fig:fas_low_08}
\end{subfigure}

\vspace{0.5em}

\begin{subfigure}[t]{0.95\textwidth}
\centering
\includegraphics[width=\linewidth]{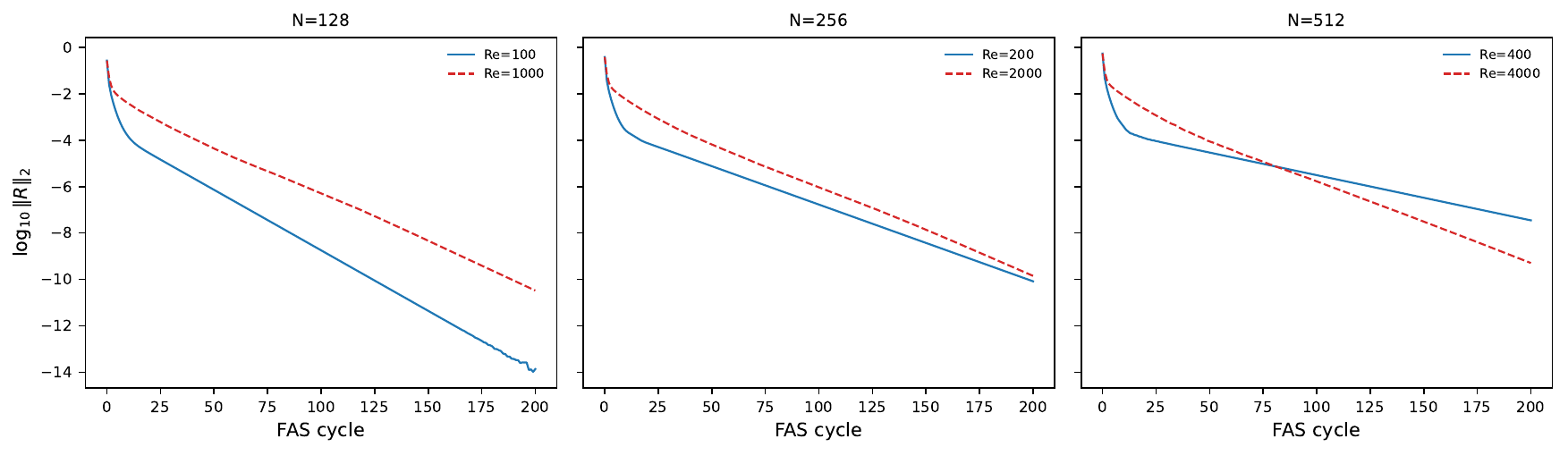}
\caption{Convergence rate depending on varying $\Rey^*$ with a fixed $N$
        }
\label{fig:fas_pairs_08}
\end{subfigure}

\caption{The evolution of the FAS-LBE residual depending on
	$N$, $\Rey^*$ (or $\omega$), with a fixed $\lambda=0.8$.
    }
\label{fig:fas_parameter_dependence_08}
\end{figure}

The effect of the smoother weight, $\lambda$, is further examined by repeating the
same set of tests
with a different value of $\lambda=0.8$ (as opposed to 0.6). The results are shown in 
\cref{fig:fas_parameter_dependence_08}.
With a larger value of $\lambda$, the overall residual
convergence behavior remains the same qualitatively, albeit some quantitative 
variations.
In contrast with the previous case of $\lambda = 0.6$, it can be seen that 
the value of $\lambda$ not only affects the absolute convergence rate of the
residual, but also its relative performance with two different values of 
$\omega$. For example, on a small mesh of $N=128$ the
lower-$\Rey^*$ case converges substantially faster with $\lambda=0.8$, 
whereas on a larger mesh of $N=512$, the higher-$\Rey^*$ case exhibits 
a faster asymptotic convergence rate.
These tests indicate that the asymptotic convergence rate of the FAS-LBE
is affected by the resolution $N$, the relaxation frequency $\omega$
(or the viscosity $\nu$), and $\lambda$. This is understandable given
the FAS-LBE involves all these parameters in its algorithmic 
construction.

Overall, the aforementioned 
results attest three features 
of the nonlinear formulation with the residual-driven LBE.
First of all, the FAS-LBE converges to the same
discrete fixed point as the primitive LBE.  Secondly, the converged
steady solution retains the accuracy of the underlying LBE discretization.
Finally, the FAS-LBE significantly enhances the computational efficiency
of the LBE algorithm when steady-state calculations are needed.
It can reduce the computational cost for steady solutions by 
orders of magnitude in the tested cases, although its asymptotic convergence 
rate remains sensitive to $\omega$ and $\lambda$.
We noted that the above tests involved only a simple Cartesian geometry, and
this motivates us to carry out more complicated tests in what follows.

\subsection{2D flow past a cylinder at $\Rey=20$}
\label{sec:dfg}

The 2D flow past a circular cylinder at $\Rey=20$ \cite{SchaeferTurek1996} 
will be used to assess the proposed FAS-LBE solver in the presence of 
a curved immersed boundary and nontrivial inflow and outflow conditions. 
In addition to 
a quantitative validation,
it may also be
used to examine the compatibility of FAS with the preconditioning strategy 
introduced in \cref{sec:solvers}.

The Reynolds number $\Rey$ is based on the cylinder diameter $D$ and the bulk
inlet velocity $\bar{U}$:
\begin{equation}
    \Rey = \frac{\bar{U} D}{\nu}=20.
\end{equation}
In what follows, we use the resolution $D = 80 \delta_x$, and 
the size of the rectangular computational domain is 
$L \times H = 22D\times4.1D = 1760 \delta_x \times 328 \delta_x$.
The cylinder center is located at $(x_c/D,\, y_c/D)=(2, 2)$.
A parabolic velocity profile is imposed
at the inlet:
\begin{equation}
    u(0,y)
    =
    4 U_m \frac{y(H-y)}{H^2},
    \qquad
    v(0,\, y)= 0 ,
\end{equation}
with the mean velocity $\bar{U}= 2 U_m / 3$. The
stationary no-slip conditions are imposed on the upper and lower 
boundaries of the domain.
The so-called do-nothing (copying) boundary conditions are imposed at the 
outlet $x=L$ 
(cf., \emph{e.g.}, \cite{junk2008outflow}).
The curved cylinder surface is treated with
the interpolated bounce-back (IBB) scheme \cite{bouzidi2001momentum}. 
The inlet Mach number is $\May=U_m/c_s = \sqrt{3} U_m$, and the viscosity is 
adjusted when $\May$ is varied to maintain a constant $\Rey$.
The hydrodynamic force acting on the cylinder is decomposed into the
streamwise and transverse components, denoted by $\FD$ and $\FL$,
respectively. They are evaluated using the momentum
exchange method \cite{bouzidi2001momentum}. 
The corresponding drag and lift coefficients are defined as
\begin{equation}
    \CD =
    \frac{\FD}
    {\frac{1}{2}\rho_0 \bar{U}^2 D},
    \qquad
    \CL =
    \frac{\FL}
    {\frac{1}{2}\rho_0 \bar{U}^2 D},
\end{equation}
where $\rho_0 = 1$ is the reference density.
The pressure difference is
evaluated between the front and rear points of the cylinder along its
horizontal centerline:
\begin{equation}
    \Delta p =
    \hat{p}\left(x_c-D/2,\, y_c \right)
    -
    \hat{p}\left(x_c+D/2,\, y_c\right).
\end{equation}

A four-level V-cycle MG scheme is employed for this case. For
the FAS-LBE, the multigrid parameters are set to
\[
\nu_1=\nu_2=3,\qquad \nu_3=40,\qquad \lambda=0.8,
\]
whereas FAS-PLBE uses
\[
\nu_1=\nu_2=2,\qquad \nu_3=20,\qquad \lambda=0.2.
\]
For the PLBE calculations, $\gamma$ is selected such that
$\May/\sqrt{\gamma}=0.2$.

\begin{figure}[htbp]
\centering
\includegraphics[width=0.75\linewidth]{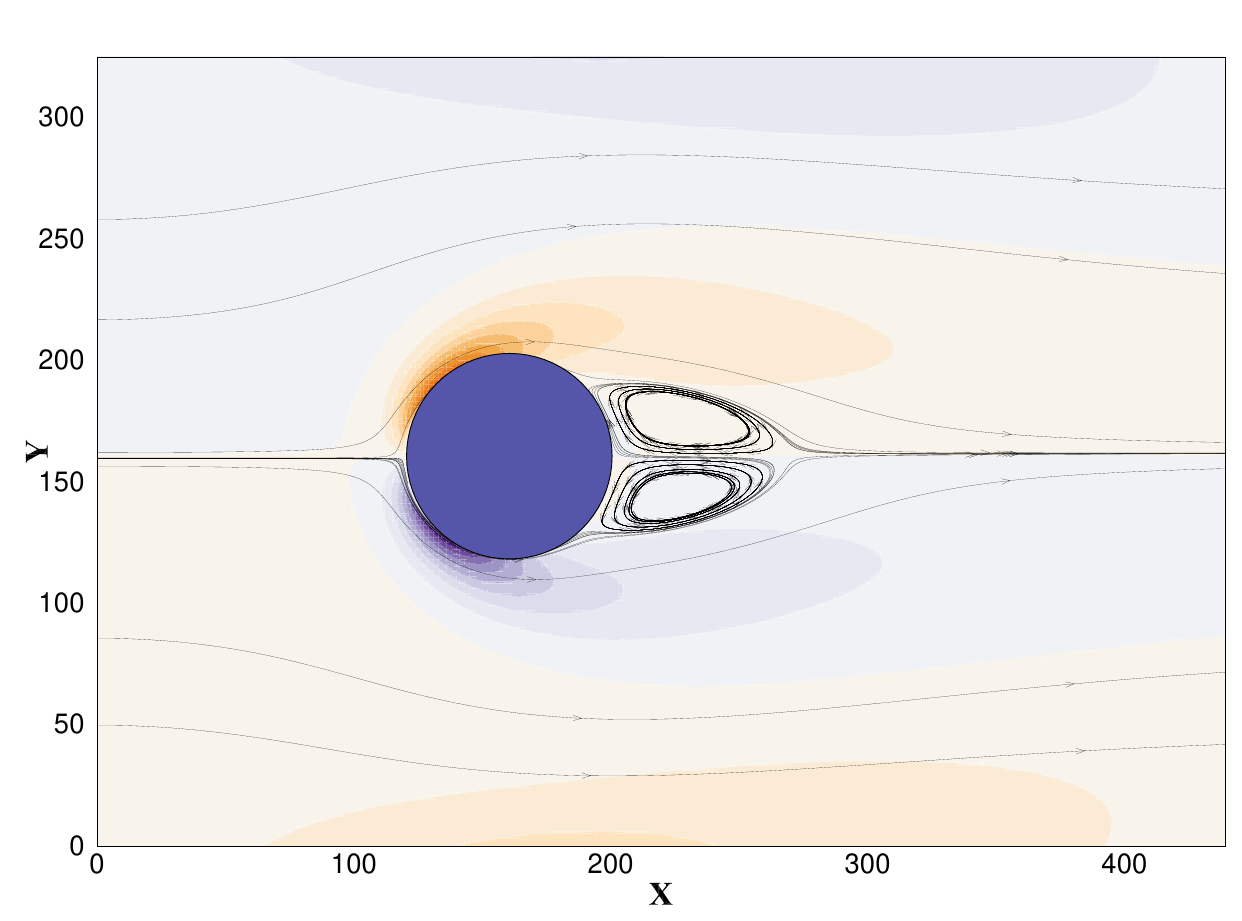}
\caption{2D flow past a cylinder at $\Rey=20$. The streamlines obtained with
the FAS-LBE scheme.}
\label{fig:dfg_flow}
\end{figure}

\cref{fig:dfg_flow} shows 
the streamlines obtained with
FAS-LBE. The characteristic steady recirculation bubbles behind
the cylinder are clearly captured. 
A more quantitative validation is to compare 
the computed drag and lift coefficients and the pressure difference, 
as compiled in \cref{tab:dfg_validation}. 
With $\May=0.1$, the $\CD$ and $\CL$ obtained with the FAS-LBE agree with
their high-order spectral counterparts \cite{Nabh1998}, and their relative differences 
in $\CD$, $\CL$, and $\Delta p$ are 0.11\%, 2.64\%, and 0.06\%, respectively. 
Our results are also consistent with the previous results obtained with 
the LBE scheme \cite{Filippova1998jcp,junk2008outflow}.
This test further confirms the fact that the FAS-LBE retains the accuracy
of the underlying LBE discretization for flows involving curved boundaries.

\begin{table}[htbp]
\centering
\caption{2D flow past a cylinder at $\Rey=20$. The 
results of $\CD$, $\CL$, and $\Delta p$. 
}
\label{tab:dfg_validation}
\begin{tabular}{lccc}
\toprule
Source & $\CD$ & $\CL$ & $\Delta p$ \\
\midrule
high-order spectral \cite{Nabh1998}
            & 5.579535 & 0.010619 & 0.117520 \\
LBGK \cite{Filippova1998jcp}
            & 5.520000 & 0.011000 & 0.116000 \\
LBGK \cite{junk2008outflow}
            & 5.582000 & 0.011100 & 0.117400 \\
FAS-LBE $\May=0.1$
            & 5.585905 & 0.010339 & 0.117588 \\
\bottomrule
\end{tabular}
\end{table}

We further examine the computational efficiency using the 2D flow past a cylinder 
at $\Rey=20$ as well as the composability of the multilevel framework. 
As discussed in \cref{sec:solvers}, the PLBE can be either employed as a stand-alone
component for iteration or incorporated into FAS as a nonlinear smoother.
Thus a FAS-PLBE scheme, as opposed to the FAS-LBE, may be easily implemented.
The FAS-LBE and the FAS-PLBE may be used to dissect
two complementary acceleration mechanisms: preconditioning at the
single-grid level and nonlinear multigrid acceleration across spatial
scales, and the outcome may help us formulate a better strategy for
enhancing efficiency.

\begin{figure}[htbp]
\centering
\includegraphics[width=0.75\linewidth]{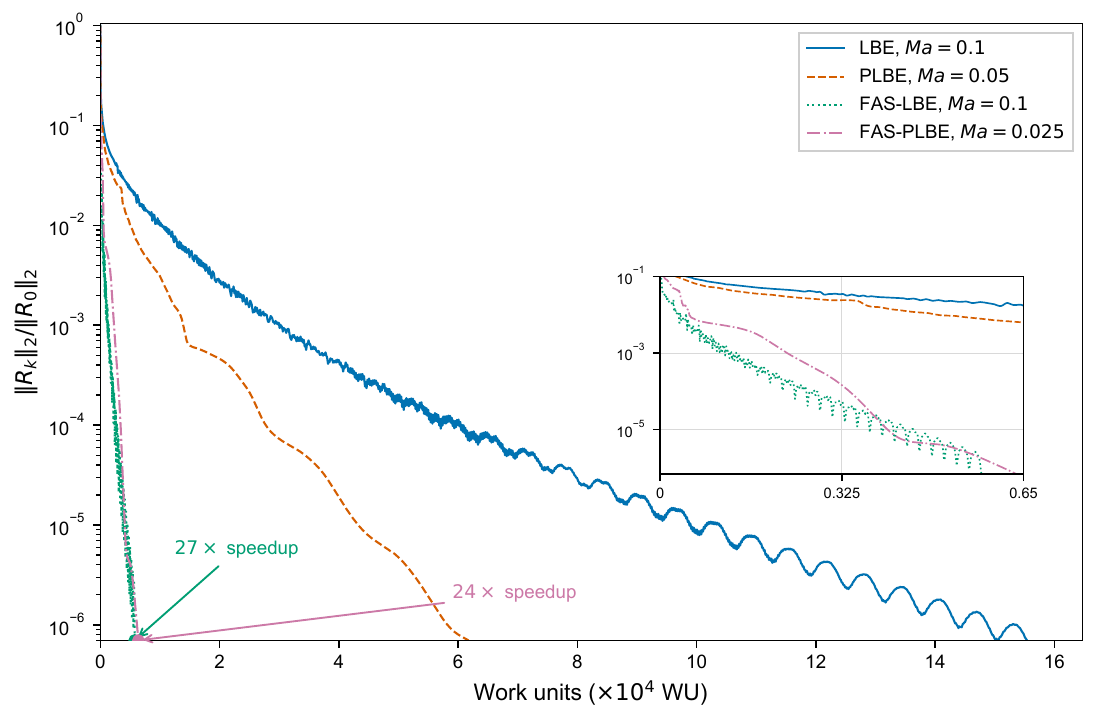}
\caption{The evolution of the residuals with the primitive LBE, PLBE, FAS-LBE, 
and FAS-PLBE for the 2D flow past a cylinder at $\Rey=20$.
}
\label{fig:dfg_efficiency}
\end{figure}

\cref{fig:dfg_efficiency} shows the evolution of the residuals in 
the primitive LBE, PLBE, FAS-LBE, and FAS-PLBE schemes vs. the work units (WUs). Clearly, both preconditioning and nonlinear multigrid 
solver accelerate computational speed tremendously, by orders of magnitude. 
For this particular test, the FAS-LBE calculation consumed 584 V-cycles, 
or equivalently 5744 WUs, while the FAS-PLBE cost 898 V-cycles or 
6364 WUs; the speedup of the FAS-LBE and FAS-PLBE is 27 and 24,
respectively, with respect to the primitive LBE iteration. One interesting
point to note is the FAS-PLBE calculation was done with $\May = 0.025$,
while the LBE was done with $\May= 0.1$. That is, the FAS-PLBE calculation
was done with the CFL number which is one quarter of that in the corresponding 
LBE calculation. The LBE calculation with $\May=0.025$ would take too long. 
With a smaller $\May$, the compressibility error is reduced. Therefore, the FAS-PLBE 
not only enhances the efficiency, but also improves the accuracy. 
We have also demonstrated the composability of the proposed residual-driven
LBE framework, that is, various nonlinear solvers may be easily wrapped around
the LBE kernel to achieve desirable enhancements. It should be stressed that
the ease of implementing the residual-driven LBE with nonlinear solvers cannot
be overstated --- this fully conforms with the simplicity of the LBE algorithm.

\subsection{2D flow past a NACA0012 airfoil at $\Rey=1000$}
\label{sec:airfoil}

The final test to be conducted is the steady external flow past a
NACA0012 airfoil with a nonzero angle of attack. This case differs from the
previous ones because it is an external flow as opposed to the internal ones,
and involves a very slender body as opposed to a bluff one. 
The Reynolds number based on the chord length $L_c$ and the freestream
velocity $U_\infty$ is
\begin{equation}
    \Rey=\frac{U_\infty L_c}{\nu}=1000,
\end{equation}
and the angle of attack is $\alpha=5^\circ$. The
computational domain size is $L \times H = 20L_c\times10L_c$. The airfoil's
leading edge is located at $5 L_c$ downstream of the inlet.
A uniform velocity condition is prescribed at the inlet, a
constant pressure condition is imposed at the outlet, and 
the free-slip conditions are used
at the upper and lower boundaries. The airfoil surface is represented
by the IBB scheme. The drag and lift coefficients are defined as
\begin{equation}
    \CD = \frac{\FD}{\frac{1}{2}\rho_0 U_\infty^2 L_c} ,
    \qquad
    \CL = \frac{\FL}{\frac{1}{2}\rho_0 U_\infty^2 L_c} .
\end{equation}
The chord is resolved by
$L_c=160\delta_x$, resulting in a domain size of $3200 \delta_x \times1600\delta_x$, and the
freestream Mach number is set to $\May_\infty = \sqrt{3} U_\infty = 0.05$.

\begin{figure}[htbp]
\centering
\includegraphics[width=0.95\linewidth]{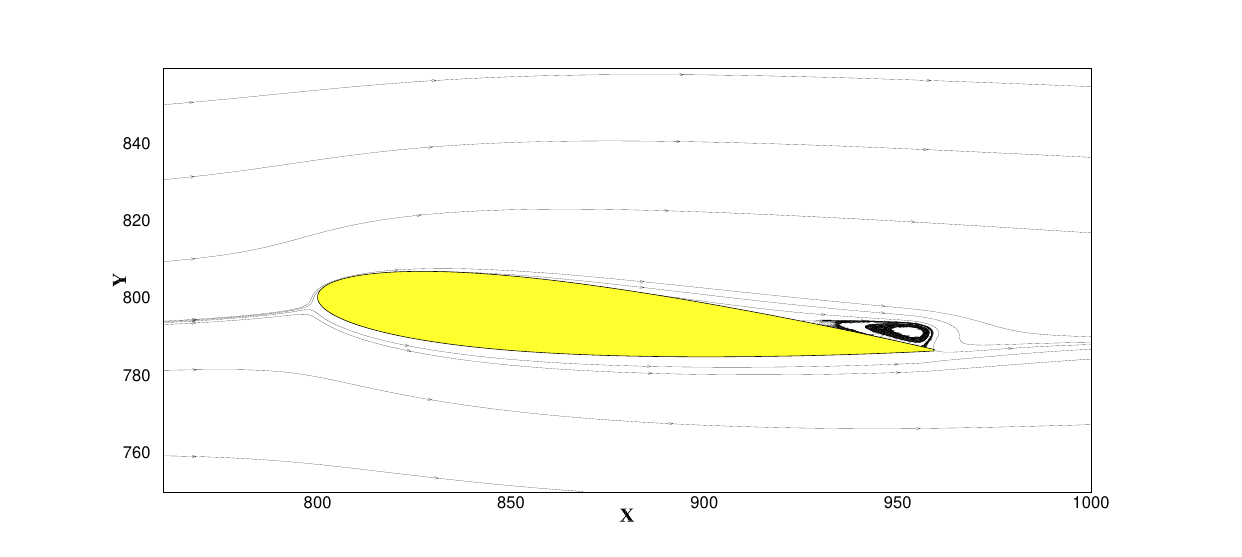}
\caption{2D flow past a NACA0012 airfoil at $\Rey=1000$, $\alpha = 5^\circ$, 
	$\May_\infty = 0.05$, and $L_c = 160 \delta_x$. 
	The streamlines around the airfoil obtained with the FAS-LBE. 
}
\label{fig:airfoil_flow}
\end{figure}

\cref{fig:airfoil_flow}
shows the converged streamline pattern around the airfoil obtained with the FAS-LBE.
The values of $\CD$ and $\CL$ are tabulated in \cref{tab:airfoil_validation} and
compared with those reported from a recent work \cite{xie2022implicit}, which used a
second-order FV method with implicit time marching scheme on cut-cell, 
non-uniform Cartesian structured mesh, and the finest grid size around the
airfoil is a rectangle $(\delta_x,\, \delta_y) = (L_c/100,\, L_c/400)$. 
The domain size in our simulation is identical to that used in 
\cite{xie2022implicit}, and the number of grid points in our calculation is
$3200 \times 1600$ on a uniform mesh, while that in \cite{xie2022implicit}
is $500 \times 300$ on a nonuniform mesh. Our results agree with the reference data in \cite{xie2022implicit} to the two significant digits reported therein.

\begin{table}[htbp]
\centering
\caption{2D flow past a NACA0012 airfoil with $\alpha=5^\circ$ at $\Rey=1000$.
	$\CD$ and $\CL$.
}
\label{tab:airfoil_validation}
\begin{tabular}{lll}
\toprule
  \multicolumn{1}{c}{Method} & 
  \multicolumn{1}{c}{$\CD$} & 
  \multicolumn{1}{c}{$\CL$} \\
\midrule
  cut-cell Cartesian 2nd-order FV \cite{xie2022implicit}
            & 0.13 &  0.25 \\
  FAS-LBE   & 0.1296  & 0.2525  \\
\bottomrule
\end{tabular}
\end{table}

The efficiencies of various FAS-LBE schemes are also to be assessed for this case.
A coarser resolution of $L_c = 80 \delta_x$ is used solely for the purpose of
assessing efficiencies, and the FAS includes a four-level FAS V-cycle. 
For the FAS-LBE and FAS-PLBE the parameters in the MG implementation are:
$$
\nu_1=\nu_2=3,\qquad \nu_3=20,\qquad \lambda=0.4.
$$
In the PLBE, $\gamma$ is chosen such that
$\May = 0.1 \sqrt{\gamma}$.

\begin{figure}[htbp]
\centering
\includegraphics[width=0.82\linewidth]{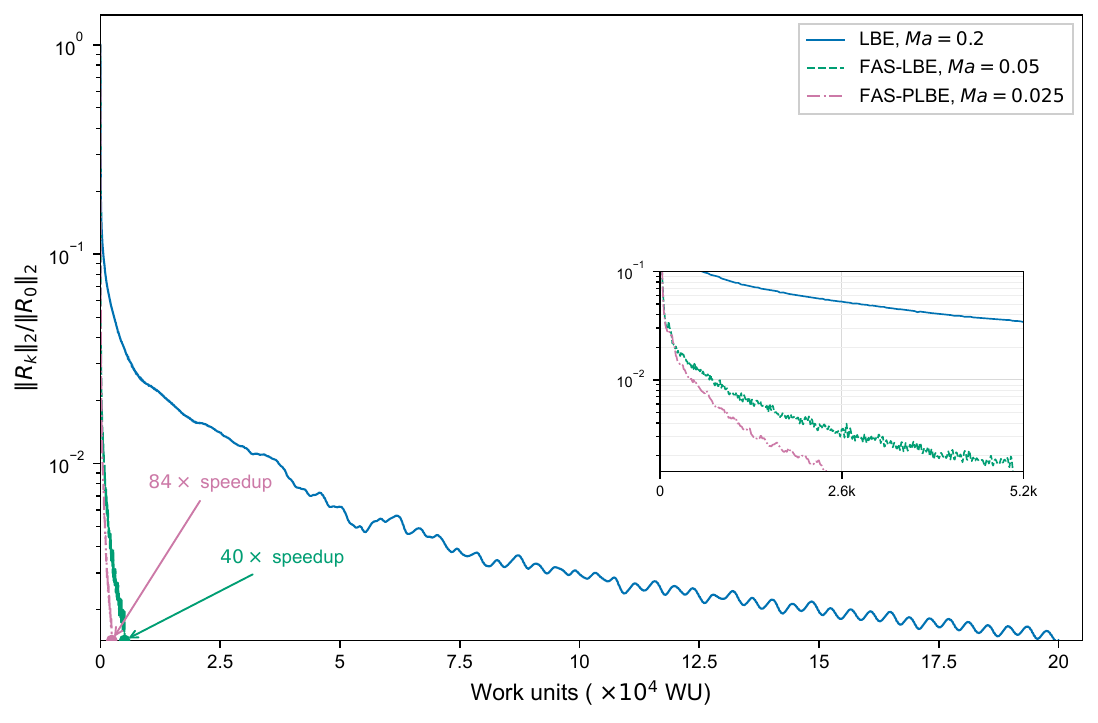}
\caption{2D flow past a NACA0012 airfoil at $\Rey=1000$ and $\alpha=5^\circ$.
The evolution of the residuals vs. WUs.
}
\label{fig:airfoil_efficiency}
\end{figure}

\cref{fig:airfoil_efficiency} shows the evolution of the residuals with
the FAS-LBE, FAS-PLBE, and primitive LBE iterations vs. WUs.
The values of the Mach number $\May$ used in the primitive LBE, 
the FAS-LBE, and the FAS-PLBE iterations are 0.2, 0.05, and 0.025, respectively.
This test clearly demonstrated the efficacy of the FAS-LBE solvers. 
In particular, the FAS-PLBE calculation took 234 V-cycles or 2382 WUs,
the FAS-LBE did 496 V-cycles or 5048 WUs, thus the speedups with 
the FAS-PLBE and the FAS-LBE are about 84 and 40 times with respect to 
the primitive LBE, respectively, in spite of the fact that the Mach 
numbers used in the FAS-PLBE and FAS-LBE calculations are much smaller than that in the 
primitive LBE --- the FAS-LBE and FAS-PLBE can effectively 
reduce the compressibility error \cite{HeXY1997jsp88,guo2004preconditioned}
and enhance efficiency simultaneously. In other words, the FAS-PLBE 
may alleviate the low-Mach-number crisis numerically in the LBE.

\section{Conclusions}
\label{sec:con}

This work developed a residual-driven LBE framework to freely incorporate
arbitrary nonlinear iterative solution methods to significantly enhance
the computational efficiency and accuracy of the LBE algorithm. The proposed framework fully segregates the nonlinear iterative solver from the internal LBE kernel, around which the solver operates; this structure thus completely
preserves the primitive LBE spatial discretizations and other key
favorable numerical characteristics. By construction, the nonlinear solver 
in this proposed framework may be freely replaced and chosen only for its 
desirable numerical features and characteristics. In this setting,
the residual computed from a complete LBE core including the collision-streaming, boundary conditions,
and possible forcing and interaction effects is used by any external 
nonlinear solver wrapping around the LBE kernel, while the LBE core algorithm serves
as an internal iterative means which is not necessarily associated with an
evolution in accord with physical-time, which may be deferred to the outer
solver of our choice. This feasibility to incorporate sophisticated numerical
techniques in the LBE algorithm has been fully displayed in the present 
work with the implementations of the nonlinear solvers based on weighted fixed-point
iteration, full approximation scheme (FAS), and Jacobian-free
Newton--Krylov methods in the proposed residual-driven LBE framework.
The composability of the proposed framework is explicitly demonstrated as well;
that is, different and desirable numerical techniques may be easily incorporated 
into the proposed framework as a stand-alone component, which is a very significant
feature.

Several numerical tests have been carried out to validate the proposed 
residual-driven LBE framework incorporated with a number of 
aforementioned nonlinear iterative solvers. The numerical tests
are all two-dimensional steady flows, including the lid-driven cavity flow at
$\Rey=1000$ and 5000; the flow past a circular cylinder at $\Rey=20$;
and the flow past a NACA0012 airfoil at $\Rey=1000$ and an angle of 
attack $\alpha=5^\circ$. These test cases consider different boundary
conditions and the results are compared with existing data to validate the
efficacy and efficiency of the proposed residual-driven LBE framework
wrapped with external nonlinear solvers. Without altering or sacrificing
the favorable characteristics of
the LBE discretization, various FAS solvers incorporated with the 
residual-driven LBE have invariably enhanced computational efficiency 
significantly, by several orders of magnitude in some cases. 
The low-Mach-number tests further demonstrate the efficacy of combining 
FAS with PLBE to reduce the computational work required for steady solutions at very low Mach numbers. 
This result suggests that the FAS-PLBE solver 
may alleviate the low-Mach-number crisis associated with the primitive LBE.
It should also be pointed out that the low-Mach-number crisis is not unique 
to the LBE; it is a common pathological issue shared by schemes based 
on artificial compressibility methodology. Our results demonstrate both the feasibility of incorporating general nonlinear iterative methods into the LBE and the substantial improvements in computational efficiency achieved by the proposed residual-driven framework.
This work further strengthens the efficacy and efficiency of the LBE as a CFD tool.

While the collision model used in the present work is the lattice BGK model,
the proposed framework is fully independent of any specific collision model in 
the LBE --- it only requires a fully encapsulated kernel with necessary components
for a complete updating step of the underlying system, which is also capable of
providing the residual to the outer nonlinear solver. Of course, the specifics of the collision model may affect, for example, the efficiency of the solver as a whole.
Secondly, the numerical tests carried out here are all steady flows in two dimensions.
However, extension of the LBE to three dimensions is a straightforward matter
\cite{dHumieres2002ptrsla, lallemand2021lattice}. In addition, a time-dependent 
solver built upon the framework can also be constructed similarly; in fact,
a preliminary feasibility study on a dual-time stepping method has been 
carried out already \cite{JiLuo2026DTSLBE}. The present work only represents
a first step to open the unbounded possibilities for incorporating advanced
numerics into the LBE method, and much is left to be done in future work,
including implicit schemes for unsteady flows in 3D and more.

\section*{Acknowledgements}

XZJ would like to acknowledge the support of the Science Challenge Project
under the Grant TZ2026004.

\appendix

\section{Hydrodynamic equations and error analysis}
\label{app:ce}

The macroscopic equations for mass and momentum conservation are to be derived from the
LBE in this section. Specifically, the physical time step $\delta_T$, the lattice 
time step $\delta_t$, and the nonlinear residual tolerance $\varepsilon_{\ssstyle \mathrm{NL}}$ 
are all to be considered.

\subsection{The tensors based on the discrete velocity set $\mathbb{X}_q$}
\label{app:lattice-identities}

Given the discrete velocity set $\mathbb{X}_q$ and the associated
weight coefficients $\{ w_\alpha \}$,
the $n$-th order tensor:
\begin{subequations}
\begin{align}
	&
	\tensorfont{E}^{(n)}
	:= 
	\sum_{\alpha}
	w_\alpha \bm{\xi_\alpha}^n
	:=
	\sum_\alpha w_\alpha \underbrace{\bm{\xi_\alpha}\cdots\bm{\xi_\alpha}}_{n} ,
	\\
	&
	E^{(n)}_{i_1 i_2 \cdots i_n}
	=
	\sum_{\alpha=0}^{8} w_\alpha
	\xi_{\alpha,i_1} \cdots \xi_{\alpha,i_n},
	\qquad i_1,\ldots,i_n\in\{1,\, 2\} = \{x,\, y\} .
\label{eq:lattice-tensor}
\end{align}
\end{subequations}
Specifically for the D2Q9 model:
\begin{subequations}
\begin{align}
	&
  E^{(0)}
  =\sum_{\alpha=0}^{8}w_\alpha=1,\\
  &
  E^{(2)}_{ij}
  =\sum_{\alpha=0}^{8}w_\alpha\xi_{\alpha,i}\xi_{\alpha,j}
  =\cs^2\delta_{ij},\\
	&
	E^{(4)}_{ijkl}
	=\sum_{\alpha=0}^{8} w_\alpha
    \xi_{\alpha,i}\xi_{\alpha,j}\xi_{\alpha,k}\xi_{\alpha,l}
  	=\cs^4 \left(
      \delta_{ij}\delta_{kl}
     +\delta_{ik}\delta_{jl}
     +\delta_{il}\delta_{jk} \right),\\
	&
	\tensorfont{E}^{(2m+1)} = \bm{0} .
\end{align}
\label{eq:lattice-isotropy}
\end{subequations}

For the equilibrium of \cref{eqn:fzero}
in \cref{sec:residual}:
\begin{subequations}
\begin{align}
	&
  {\textstyle \sum_{\alpha=0}^{8}} \fazero
  =\frac{\hat p}{\cs^2} ,
  \qquad
  {\textstyle \sum_{\alpha=0}^{8}} \bm{\xi}_\alpha \fazero
  =\bm{u},
  \\
  &
  {\textstyle \sum_{\alpha=0}^{8}} \xi_{\alpha,i}\xi_{\alpha,j} \fazero
  =
  \hat p\, \delta_{ij}+u_i u_j ,
\label{eq:eq-second-moment}\\
  &
  {\textstyle \sum_{\alpha=0}^{8}}
  \xi_{\alpha,i}\xi_{\alpha,j}\xi_{\alpha,k} \fazero
  = 
  \cs^2\left(
      \delta_{ij}u_k+\delta_{ik}u_j+\delta_{jk}u_i\right).
\label{eq:eq-third-moment}
\end{align}
\end{subequations}

\subsection{Asymptotic expansion}
\label{app:ce-expansion}

Denote $\fceva(\bm{r_j},\, t) := f_\alpha( \bm{r_j} - \bm{\xi}_\alpha \delta_t, \, t)$
the distribution pulled from the upwind lattice node
$\bm{r_j}-\bm{\xi}_\alpha\delta_t$ along $-\bm{\xi}_\alpha$.
Assume $\kappa$ to be the order of convergence for the physical-time discretization:
\begin{equation}
  \tensorfont{D}_{T}(f_\alpha)
  =
  \frac{\diff f_\alpha}{\diff t}
   + O(\delta_T^\kappa),
\end{equation}
so that the nonlinear equation has a residual of order
$O(\varepsilon_{\ssstyle \mathrm{NL}})$ in the specified normalization.  
Taylor expansion of the pull representation of \cref{eq:general-residual} gives
\begin{equation}
\begin{aligned}
  \frac{\partial \fceva}{\partial t}
   =
   &
   \Omegaceva -
   \bm{\xi}_\alpha \!\cdot\! \bm{\nabla}
  \left[ 1 + \delta_t \left( {\textstyle \frac{1}{2}} 
  \bm{\xi}_\alpha \!\cdot\! \bm{\nabla} + \partial_t 
  \right) \right] \fceva
   \\
   & + O(\delta_t^2+\delta_T^\kappa +\varepsilon_{\ssstyle \mathrm{NL}}).
\end{aligned}
\label{eq:ce-taylor}
\end{equation}
Expand the distribution as
\begin{equation}
  f_\alpha
  = \fazero
   +\delta_t \faone
   +\delta_t^2 \fatwo+
   \cdots ,
\label{eq:ce-distribution}
\end{equation}
and define the physical-time transport operator
\begin{equation}
  \tensorfont{K}_\alpha
  :=
  \frac{\diff}{\diff t} +
  \bm{\xi}_\alpha \!\cdot\! \bm{\nabla} .
\end{equation}
Substitution of \cref{eq:ce-distribution} into \cref{eq:ce-taylor}, 
followed by removal of the common pull shift, yields
\begin{equation}
\begin{aligned}
  &
  \tensorfont{K}_\alpha \fazero
  + \delta_t
  \left[
      \tensorfont{K}_\alpha \faone
     +\frac{1}{2}(\bm{\xi}_\alpha\cdot \bm{\nabla} )^2 \fazero
  +  \bm{\xi}_\alpha \!\cdot\! \bm{\nabla}
       \frac{\partial \fazero}{\partial t}
	\right]
	\\
  = & 
  - \omega \faone
   + O(\delta_t^2+\delta_T^\kappa+\varepsilon_{\ssstyle \mathrm{NL}}).
\end{aligned}
\label{eq:ce-first-form}
\end{equation}
At leading order,
\begin{equation}
  \faone
  =-\frac{1}{\omega} \tensorfont{K}_\alpha \fazero
   +O(\delta_t+\delta_T^\kappa+\varepsilon_{\ssstyle \mathrm{NL}}).
\label{eq:ce-first-order-noneq}
\end{equation}
Substituting \cref{eq:ce-first-order-noneq} in \cref{eq:ce-first-form}
leads to the following equality used to derive the moment equations later:
\begin{equation}
\begin{aligned}
  \tensorfont{K}_\alpha \fazero
  =
  & -\omega \faone
  + \frac{\nu}{\cs^2}
    \left[
      ( \bm{\xi}_\alpha \!\cdot\! \bm{\nabla} )^2 
      +2\bm{\xi}_\alpha \!\cdot\! \bm{\nabla} \partial_t
    \right] \fazero
	\\
  &+\frac{\delta_t}{\omega}
    \frac{\partial^2 \fazero}{\partial t^2}
   +O(\delta_t^2+\delta_T^\kappa+\varepsilon_{\ssstyle \mathrm{NL}}).
\end{aligned}
\label{eq:ce-final-form}
\end{equation}

\subsection{Continuity equation}\label{app:ce-mass}

Assume the incompressible scaling such that:
\begin{equation}
  \bm{u} = O(\May),
  \qquad
  \hat p=O(\May^2),
\end{equation}
and the compatibility conditions:
\begin{equation}
  \sum_{\alpha=0}^{8}f_\alpha^{(m)} = 0,
  \qquad
  \sum_{\alpha=0}^{8}\qvector{\xi}_\alpha f_\alpha^{(m)}
  = \bm{0} ,
  \qquad m \geq 1 .
\label{eq:ce-noneq-moments}
\end{equation}
Taking the zeroth moment of \cref{eq:ce-final-form} gives
\begin{equation}
  \frac{1}{\cs^2} \frac{\partial \hat{p}}{\partial t}
  +\bm{\nabla} \!\cdot\! \bm{u}
  = \frac{\nu}{\cs^2}P
   +\frac{\delta_t}{\cs^2\omega}
      \frac{\partial^2 \hat{p}}{\partial t^2}
   +O(\delta_t^2+\delta_T^\kappa+\varepsilon_{\ssstyle \mathrm{NL}}),
\label{eq:ce-mass}
\end{equation}
where
\begin{equation}
\begin{aligned}
  P
  &:=\partial_i \partial_j
       \sum_{\alpha=0}^{8}\xi_{\alpha,i}\xi_{\alpha,j} \fazero
    +2 \frac{\partial}{\partial t}
       \bm\nabla\cdot
       \sum_{\alpha=0}^{8} \bm{\xi}_\alpha \fazero
	   \\
  &=\bm{\nabla} \!\cdot\! 
  \left[ \bm{\nabla} \!\cdot\! (\bm{u u}) \right]
    +\Delta\hat p
    +2\frac{\partial}{\partial t}(\bm{\nabla} \!\cdot\! \bm{u} )
   =O(\May^2).
\end{aligned}
\label{eq:ce-P}
\end{equation}
Consequently,
\begin{equation}
  \bm{\nabla} \!\cdot\! \bm{u}
  = O(\May^2\delta_t+\delta_t^2
    +\delta_T^\kappa+\varepsilon_{\ssstyle \mathrm{NL}}).
\label{eq:ce-divergence}
\end{equation}

\subsection{Momentum equation}\label{app:ce-momentum}

Taking the first moment of \cref{eq:ce-final-form} and using
\cref{eq:eq-second-moment,eq:eq-third-moment} gives
\begin{equation}
\begin{aligned}
  \frac{\partial \bm{u}}{\partial t}
  +\bm{\nabla} \!\cdot\! (\bm{u u})
  +\bm{\nabla} \hat{p}
  = \,
  \frac{\nu}{\cs^2}
  \bm{Q}
   +\frac{\delta_t}{\omega}
     \frac{\partial^2\bm{u}}{\partial t^2}
  +O(\delta_t^2+\delta_T^\kappa+\varepsilon_{\ssstyle \mathrm{NL}}),
\end{aligned}
\label{eq:ce-momentum}
\end{equation}
where
\begin{equation}
\begin{aligned}
  \bm{Q}
  &:=\cs^2\left[
       \Delta \bm{u}+2 \bm{\nabla} (\bm{\nabla} \!\cdot\! \bm{u} )
     \right]
    +2 \bm{\nabla} \!\cdot\!
      \frac{\partial}{\partial t}
      \left( \hat{p} \tensorfont{I}+\bm{u u} \right)
	  \\
  &= \cs^2\Delta \bm{u}
    + 2 \bm{\nabla} \!\cdot\!
      \frac{\partial}{\partial t}
      \left(\hat{p} \tensorfont{I}+\bm{u u} \right)
    +O(\May^2).
\end{aligned}
\label{eq:ce-Q}
\end{equation}
Let
$\tensorfont{\Pi}^{\ssstyle (0)} := \hat{p} \tensorfont{I}+\bm{u u}$
and
$\textstyle \nu = \cs^2 \delta_t (\frac{1}{\omega} -\frac{1}{2} )$.
%
%
\cref{eq:ce-momentum} can be rewritten as
\begin{equation}
\begin{aligned}
  \frac{\partial \bm{u}}{\partial t}
  +\bm{\nabla} \!\cdot\! \tensorfont{\Pi}^{\ssstyle (0)}
  = \,
  & \nu \Delta \bm{u}
  + 
   \frac{\delta_t}{\omega} \frac{\partial^2 \bm{u}}{\partial t^2}
   + \frac{2 \nu}{\cs^2} 
   \frac{\partial}{\partial t} 
   \bm{\nabla} \!\cdot\! \tensorfont{\Pi}^{\ssstyle (0)}
	\\
	&
   +O(\May^2\delta_t+\delta_t^2
     +\delta_T^\kappa+\varepsilon_{\ssstyle \mathrm{NL}})\\
  ={}&\nu\Delta \bm{u}
   +O(\May^2\delta_t+\delta_t^2
     +\delta_T^\kappa+\varepsilon_{\ssstyle \mathrm{NL}}).
\end{aligned}
\label{eq:ce-NS}
\end{equation}
Thus the residual-driven physical-time equation recovers
\begin{equation}
\begin{aligned}
  \bm{\nabla} \!\cdot\! \bm{u}
  &= O(\May^2\delta_t+\delta_t^2
     +\delta_T^\kappa+\varepsilon_{\ssstyle \mathrm{NL}}),\\
  \frac{\partial\bm{u}}{\partial t}
  +\bm\nabla \!\cdot\! (\bm{u u})
  &=-\bm{\nabla} \hat{p}+\nu \Delta \bm{u}
   +O(\May^2\delta_t+\delta_t^2
     +\delta_T^\kappa+\varepsilon_{\ssstyle \mathrm{NL}}).
\end{aligned}
\label{eq:ce-final-macroscopic}
\end{equation}

\subsection{The modified-equations}
\label{app:ce-comparison}

The modified equations 
may be 
obtained by the Taylor expansion of the primitive LBE in $\delta_t$ without assuming any scaling relationship
\cite{junk2005asymptotic, Dubois2008cma, Bellotti2024jcp}. 
The resulting equations for the mass and momentum conservation equations are:
\begin{equation}
  \frac{1}{\cs^2}\frac{\partial \hat{p} }{\partial t}
  +\bm{\nabla} \!\cdot\! \bm{u}
  =\frac{\nu}{\cs^2}P
   +\frac{\nu}{\cs^4}
    \frac{\partial^2 \hat{p} }{\partial t^2}
   +O(\delta_t^2),
\label{eq:ce-mass-lbe}
\end{equation}
and
\begin{equation}
\begin{aligned}
  \frac{\partial \bm{u}}{\partial t}
  +\bm{\nabla} \!\cdot\! \tensorfont{\Pi}^{\ssstyle (0)}
  = 
  \, &\nu \Delta \bm{u} 
  + \frac{\nu}{\cs^2}
  \left(
    \frac{\partial^2 \bm{u}}{\partial t^2}
  + 
    2 \frac{\partial}{\partial t}
     \bm{\nabla} \!\cdot\! \tensorfont{\Pi}^{\ssstyle (0)}
	\right)
	 \\
  &+O(\May^2\delta_t+\delta_t^2).
\end{aligned}
\label{eq:ce-momentum-lbe}
\end{equation}
\cref{eq:ce-mass} with \cref{eq:ce-mass-lbe}, and
\cref{eq:ce-NS} with \cref{eq:ce-momentum-lbe}, 
contain similar leading terms. However, the coefficients 
of higher-order temporal terms, \emph{i.e.}, $\partial^2 \hat{p}/\partial t^2$ and
$\partial^2 \bm{u}/\partial t^2$,
are not the same. Under the assumed low-Mach-number scaling, the difference enters the truncation structure beyond the leading incompressible Navier--Stokes equations.
The physical-time discretization term $O(\delta_T^\kappa)$ and the
nonlinear-residual contribution $O(\varepsilon_{\ssstyle \mathrm{NL}})$
are additional perturbations in the macroscopic equations induced by
the outer time discretization and the inexact nonlinear solve, respectively.
The latter is not automatically a bound on the solution error.

\bibliographystyle{elsarticle-num}
\bibliography{journalnames,mgbib,latticeBoltzmann}

\end{document}